\documentclass[11pt]{article}

\title{Bayesian Multiscale Finite Element Methods. Modeling missing subgrid information probabilistically }

\author{
Y. Efendiev \thanks{Department of Mathematics \& Institute for Scientific Computation (ISC),
Texas A\&M University,
College Station, Texas, USA. Email: {\tt efendiev@math.tamu.edu}.}
\and
W.T. Leung \thanks{Department of Mathematics, Texas A\&M University, College Station, TX 77843}
\and
S. W. Cheung \thanks{Department of Mathematics, Texas A\&M University, College Station, TX 77843}
\and
N. Guha\thanks{Department of Statistics, Texas A\&M University, College Station, TX 77843 \&  Institute for Scientific Computation, Texas A\&M University, College Station, TX 77843-3368}
\and
V. H. Hoang\thanks{School of Physical \& Mathematical Sciences, Nanyang Technological University, Singapore}
\and
B. Mallick\thanks{Department of Statistics, Texas A\&M University, College Station, TX 77843}
}

\usepackage{amsmath}
\usepackage{amsfonts}
\usepackage{amssymb}
\usepackage{graphicx}
\usepackage{algorithm}
\usepackage{algpseudocode}
\usepackage{mathrsfs}
\usepackage{caption}
\usepackage{subcaption}
\usepackage{pgffor}
\usepackage{amsthm}
\usepackage{placeins}
\usepackage{enumerate}
\usepackage{color}
\usepackage{bm}
\usepackage{epstopdf}

\usepackage[normalem]{ulem}

\newtheorem{remark}{\indent\sc Remark}

\newcommand{\note}[1]{{\color{black}#1}}

\begin{document}

\bibliographystyle{unsrt}

\maketitle

\begin{abstract}

In this paper, we develop a Bayesian multiscale approach based
on a multiscale finite element method. Because of scale disparity
in many multiscale applications,
computational models can not resolve all scales. Various subgrid
models are proposed to represent un-resolved scales. Here, we consider
a probabilistic approach for modeling un-resolved scales using the
Multiscale Finite Element Method
(cf., \cite{chkrebtii2016,mallick2016}).
  By representing dominant modes
using the Generalized Multiscale Finite Element, we propose
a Bayesian framework, which provides
multiple inexpensive (computable) solutions for a deterministic problem.
These approximate probabilistic solutions may not be very close to
the exact solutions and, thus, many realizations are needed.
In this way, we obtain a rigorous probabilistic description
of approximate solutions.
 In the paper, we consider
parabolic and wave equations in heterogeneous media.
In each time interval, the domain is divided into subregions.
Using residual information, we design appropriate prior and posterior
distributions. The likelihood consists of the residual minimization.
To sample from the resulting posterior distribution, we consider
several sampling strategies.
The sampling involves identifying important regions and important
degrees of freedom beyond {\it permanent} basis functions,
which are used in residual computation.
Numerical results are presented. We consider two sampling algorithms.
The first algorithm uses sequential sampling and is inexpensive.
In the second algorithm, we perform full sampling using the Gibbs sampling
algorithm, which is more accurate compared to the sequential sampling.
The main novel ingredients of our approach consist of:
defining appropriate permanent basis functions and the corresponding
residual; setting up a proper posterior distribution; and
sampling the posteriors.

\end{abstract}

\section{Introduction}

Many problems in application domains have multiple scales.
The scales in space and time are dominant in these applications
that arise in porous media, material sciences and so on.
Detailed descriptions at the finest scales often
include uncertainties due to missing information.
Moreover, there is often limited information about the solution available.
For this reason, it is desirable to compute solutions
within a probabilistic setup and estimate
associated uncertainties, which is an objective of this paper.

One of the challenges in the computation of multiscale
problems is the resolution of finest scales.
Due to the grid resolution, one can not afford many degrees of freedom in
each computational grid. Subgrid information are often modeled
stochastically even for deterministic problems. For example,
a typical approach for modeling subgrid information uses
Representative Volume Element (RVE) to compute macroscopic
quantities. However, because of uncertainties in RVE sizes and boundary
conditions, the macroscopic parameters can not be modeled deterministically.
It is advantageous in this and other multiscale applications to use
probabilistic approaches to compute the solution.
In this paper, our goal is to propose a novel modeling approach
for missing subgrid information by setting up a Bayesian formulation.

The reduced-order modeling approaches have been commonly used
in solving multiscale problems. These approaches include
homogenization and numerical homogenization methods
\cite{weh02,dur91,fish_book,eh09,fan2010adaptive, fish2012homogenization,li2008generalized,brown2013efficient}, multiscale  methods \cite{fish_book,ee03,abe07,oh05,hw97,eh09, ehw99, egw10,fish2012staggered,chung2016adaptive,efendievsparse, franca2005multiscale,nouy2004multiscale,calo2016multiscale}, and so on.
In homogenization
and numerical homogenization approaches, the macroscopic information
is computed using RVE simulations. In these approaches, the sizes of RVEs and
limited number of local information can be insufficient to compute
the solution accurately. In multiscale finite element methods,
in particularly in Generalized Multiscale Finite Element Method
(GMsFEM), multiscale basis functions are computed
to systematically
take into account missing subgrid information. In these approaches,
the missing subgrid information is represented in the form
of local multiscale basis functions. Using a few initial (dominant)
basis functions, the error can be reduced substantially.
The multiscale basis functions are computed under some assumptions.
Our objective is to propose a Bayesian formulation, which can allow
computing the multiscale solution and associated uncertainties
with a few basis functions and stochastically representing
the missing information.

The probabilistic approaches are important for problems
when one has limited information about the solution.
For example, it is common to measure the solution or averages
at some locations with some
precisions. In this case, to impose the constraint on the solution
at some time instants, one can easily use
additional constraints
in the
Bayesian framework by including additional multiplicative factors.
Furthermore,
one can include uncertainties in the media properties in the
Bayesian framework and compute the solution and the uncertainties associated
with the solution and the variations of the field parameters.
Bayesian approaches for forward and inverse problems have been
developed in many previous papers (e.g., \cite{bilionis2013multi,bilionis2013solution,marzouk2009stochastic,arnst2010identification,stuart2010inverse}).
Our approach shares some similarities with \cite{owhadi2015bayesian},
though
in our paper, we seek multiscale basis functions in a given set of basis
functions.

Our approach starts with the Generalized Multiscale Finite Element framework.
The GMsFEM was first presented in \cite{egh12}
and later investigated in several other papers
(e.g., \cite{galvis2015generalized,Ensemble, eglmsMSDG, eglp13oversampling,calo2014multiscale,chung2014adaptive,chung2015generalizedperforated,chung2015residual,chung2015online}).
It is a generalization of the MsFEM and
defines appropriate local snapshots and local spectral decompositions.
This approach adds local degrees of freedom as needed and provides
{\it numerical} macroscopic
equations for problems without scale separation and
identifies important features for multiscale problems.
Because of the local nature of proposed multiscale model reduction,
the degrees of freedom can be added adaptively based on
error estimators.
However, due to the computational
cost, one often uses
fewer basis functions. This can result to discretization
errors, which we would like to
represent in a Bayesian framework.

\begin{figure}[t]
\centering
\includegraphics[width=4in, height =3in]{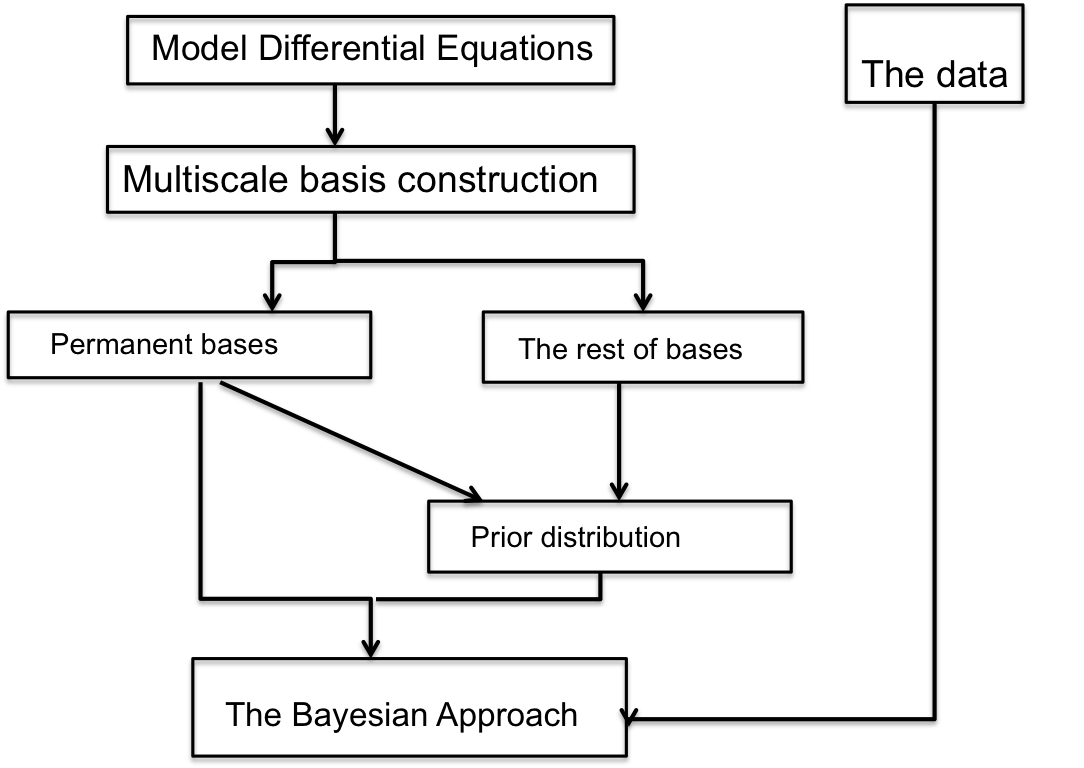}
\caption{Outline of the algorithm.}
\label{fig:mmr1}
\end{figure}

We consider the time-dependent
 equations
\[
{\partial u \over \partial t}  = \mathcal{L}(\kappa(x,t), u, \nabla u),
\]
where $\kappa(x,t)$ is a heterogeneous space-time function
and $\mathcal{L}$ is a differential operator.
Our approach starts with constructing
multiscale basis functions
and uses a few basis functions as permanent basis functions
(see Figure \ref{fig:mmr1}). It is known that these basis functions
can provide a solution approximation. Additional basis functions
are selected stochastically using the residual information.
These basis functions are selected conditioned to the selection
of subregions, which is based on the distribution of local residuals.
The latter is used as a prior distribution.

To setup the posterior distribution, we start with an approximate
solution computed with a few permanent basis functions
and compute the corresponding residual. Using the residual
information at the current time,
we define a prior distributions for the basis selection.
The likelihood includes the residual, which is minimized.
We discuss several choices for the posterior distributions
and two sampling algorithms.
The first sampling algorithm is a
sequential sampling and uses the prior distributions based
on the residual to select the realizations of the solution.
The second sampling algorithm,
full sampling,
seeks basis functions and the solution that can provided a desired error
distribution.
We note that a general flexibility of our proposed framework
allows implementing various solution strategies and
incorporating the data.
In our approach, we do not seek `''very close'' approximations
of the solutions; but rather look for many approximate solutions.
Below, we summarize our algorithm (see Figure \ref{fig:mmr1}):
\begin{itemize}

\item  compute multiscale basis functions and identify permanent basis
functions and ``the rest'' of basis functions;

\item use the residual to compute prior distribution for ``the rest'' basis
functions;

\item  setup a posterior, which
includes the residual minimization and the data;

\item  sample the posterior
distribution.

\end{itemize}

We present numerical results for parabolic and wave equations.
Our numerical results study
probabilistic approximations of the solution.
We show that the full sampling provides better
accuracy at a higher computational cost.
The number of multiscale basis functions in full sampling
depends on the threshold in the posterior and this is studied
in the paper.
The main novel ingredients of our approach are
(1) defining appropriate permanent basis functions and corresponding
residuals for priors; (2) setting up a proper posterior distribution,
the likelihood, and the prior distribution; and (3)
sampling methods, which explore these posterior distributions.

The paper is organized as follows. In Section
\ref{sec:prelim}, we present a general problem setup and
motivation. Section \ref{sec:bayes} is devoted to
our Bayesian algorithm. In this section, we describe some basic
ingredients of our algorithm and setup the posterior distribution.
In Section \ref{sec:numerics}, we present numerical
results for parabolic and wave equations.


\section{Problem setup, preliminaries, and motivation}
\label{sec:prelim}

We consider the forward model
\begin{equation}
\label{eq:main}
{\partial u \over\partial t} = \mathcal{L}(\kappa(x,t),u, \nabla u),
\end{equation}
where $\mathcal{L}$ is a multiscale differential operator.
For example, $\mathcal{L}(\kappa(x),u, \nabla u)= \text{div} (\kappa(x,t)\nabla u)$
or higher order differential operator, where
$\kappa(x,t)$ is highly oscillatory coefficients.

\begin{figure}[tb]
  \centering
  \includegraphics[scale=0.6]{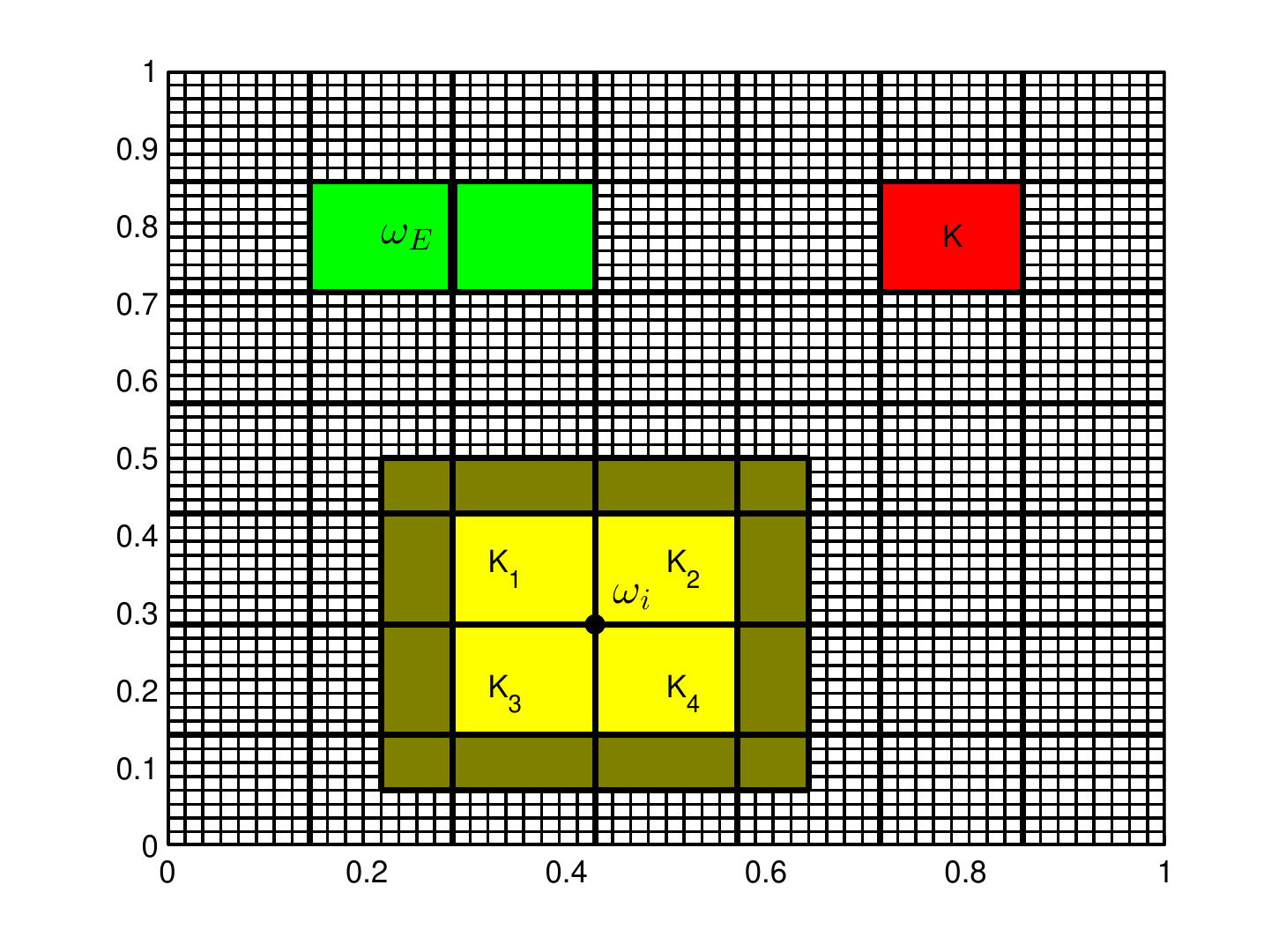}
  \caption{Illustration of fine grid, coarse grid, coarse neighborhood and oversampled domain.}
  \label{schematic_intro}
\end{figure}

We consider solving (\ref{eq:main}) on a coarse grid
(see Figure \ref{schematic_intro}
for the illustration).
Previous approaches  construct multiscale basis functions
on the coarse grid block to represent important degrees
of freedom. However, the missing degrees of freedom are not modeled
and the construction of basis functions often rely on some
assumptions. In this paper, we use a Bayesian framework and develop
a novel multiscale approach.

\subsection{Computational grid. Description of coarse and fine grids}

We introduce the notation for the coarse and fine grid.
Computations are done on a coarse grid, where fewer
basis functions are used.
The computational domain is denoted by $\Omega$, which
is partitioned by a coarse grid $\mathcal{T}^H$.
The coarse grids contain multiscale features of the problem
and require many degrees of freedom for modeling.
We denote by $N_c$, the number of nodes in the coarse grid and
by $N_e$ be the number of coarse edges.
$K$ is a generic coarse element in
$\mathcal{T}^H$.
Multiscale basis functions are computed on
a refinement of $\mathcal{T}^H$,
called a fine grid $\mathcal{T}^h$,
with mesh size $h>0$. The fine grid resolves
multiscale features of the problem.

\subsection{General idea of multiscale basis construction}

We discuss some general ideas of multiscale
basis construction, which allows building basis functions
to approximate
some important features of the solution.
We will focus on Generalized Multiscale Finite Element Method (GMsFEM).
First, we discuss a variational formulation
and then discuss the basis constructions. Some of the technical
details about basis constriction are presented in Appendix \ref{sec:app1}.

\subsection{Variational formulation}

The variational form reduces the PDE system into a system of
linear equation. The dimension of the system is proportional to the
number of basis. For very high dimensional system, as those
considered in the paper, this
leads to very large computational problems.
We
propose a Bayesian basis
selection technique.
This approach uses smaller number of basis on space-time grids
in capturing
the unresolved scale through stochastic approach.
Typical subgrid basis functions representing the
solutions over computational grids can not include all fine-grid information
of the solution space. Some
important subgrid information can be taken into account,
while  un-resolved scales and information can be
modeled in a probabilistic fashion, as proposed in the paper.

As an example, we consider the  parabolic differential equation
(\ref{eq:main}) in a space-time domain  $\Omega \times (0,T)$
and assume $u=0$ on $\partial\Omega\times(0,T)$ and
$u(x,0)=\beta(x)$ in $\Omega$. We assume the source term is $f(x)$.
We compute the solution $u_{H}$ in the
time interval $(0,T)$.
The solution
space is denoted by
$V_{H,\text{off}}^{(0,T)}$, which is a direct sum of the spaces only
containing the functions
defined on one single coarse
time interval $(T_{n-1},T_{n})$, we can decompose
the problem into a sequence of problems
and find the solution $u_{H}$ in each time interval,
denoted by $u_H^n$.
The coarse space in each
$(T_{n-1},T_{n})$
 is constructed and denoted by
$V_{H,\text{off}}^{(0,T)}=\oplus_{n=1}^{N}V_{H,\text{off}}^{(T_{n-1},T_n)}$,
 where $V_{H,\text{off}}^{(T_{n-1},T_n)}$
contains the functions having zero values
in the time interval $(0,T)$ except $(T_{n-1},T_{n})$, namely
$\forall v\in V_{H,\text{off}}^{(T_{n-1},T_n)},$
$v(\cdot,t)=0\text{ for }t\in(0,T)\backslash(T_{n-1},T_{n})$.

The coarse-grid equation is to find
\begin{equation}
  \note{
  \begin{split}
    R^n_v(u_{H}^{n}(x,t),u_{H}^{n-1}(x,t))=\\
    \int_{T_{n-1}}^{T_{n}}\int_{\Omega}\cfrac{\partial u_{H}^{n}}{\partial t}v+A(u_{H}^{n},u_H^{n-1},v)=0,
\end{split}
}
  \end{equation}
where $A$ corresponds to the PDE and
its discretization and $v$ is the test functions.
We will also investigate wave equations described in
Section \ref{sec:numerics}.
In the paper, we will consider spatial formulations, but the method
can be used for space-time method.

{\bf Example.}
In a continuous Galerkin formulation for parabolic equations
with heterogeneous coefficients $\kappa$, we seek
$u_{H}^{n}\in V_{H,\text{off}}^{(T_{n-1},T_n)}$
(where $V_{H,\text{off}}^{(T_{n-1},T_n)}$ will be defined later)
satisfying
\begin{equation}
\label{eq:space-time-FEM-coarse-decoupled}
\begin{split}
\int_{T_{n-1}}^{T_{n}}\int_{\Omega}\cfrac{\partial u_{H}^{n}}{\partial t}v+\int_{T_{n-1}}^{T_{n}}\int_{\Omega}\kappa\nabla u_{H}^{n}\cdot\nabla v+\int_{\Omega}u_{H}^{n}(x,T_{n-1}^{+})v(x,T_{n-1}^{+})\\
=  \int_{T_{n-1}}^{T_{n}}\int_{\Omega}fv+\int_{\Omega}g_{H}^{n}(x)v(x,T_{n-1}^{+}),
\end{split}
\end{equation}
for all $v\in V_{H,\text{off}}^{(T_{n-1},T_n)}$,
where
$g_{H}^{n}(\cdot)=\{u_{H}^{n-1}(\cdot,T_{n-1}^{-}) \text{ for }n\geq1;
\beta(\cdot)  \text{ for }n=0\}$,
and $F(\alpha^+)$ and $F(\alpha^-)$ denote the right hand and left hand limits of $F$ at $\alpha$ respectively.
Then, the solution $u_{H}$ of the problem in $(0,T)$
is the direct sum of all these $u_{H}^{n}$'s, that is $u_{H}=\oplus_{n=1}^{N}u_{H}^{n}$.
In multiscale simulations, our
objective is to define multiscale basis functions and minimize
$R^n$ in the space of these basis functions. In our current approach,
we would like to setup a Bayesian framework and sample a probability
distribution related to $R^n$.

\begin{remark}
We note that the residual can be written in a discrete form.
In general, the discrete form is
\[
R( u^{n+1}_{\text disc},  u^{n}_{\text disc})=\Psi_v M_{\text fine}  \Phi_u u^{n+1}_{\text disc} - \Psi_v M_{\text fine}  \Phi_u u^{n}_{\text disc}+
\Delta t \Psi_v A_{\text fine} \Phi_u u^{n+1}_{\text disc},
\]
where $M_{\text fine}$ and $ A_{\text fine}$ are fine-grid mass and stiffness
matrices, $\Psi_v$ consists of
 the test space basis vectors, $\Phi_u$ consists of
 the trial basis vectors,
and $ u^{n}_{\text disc}$ is a discretized solution at $n$th time step.
The test functions are important for the stability in the proposed method.
In the paper, we will mostly use the test spaces that correspond
to snapshot spaces, but in general, one can use a test space consisting
of all fine-grid functions.
\end{remark}

Next, we discuss multiscale basis function construction procedure
with details described in Appendix \ref{sec:app1}.

\subsection{Multiscale basis functions and snapshot spaces.}

The main idea of a multiscale approach is to systematically select
important degrees of freedom for the solution in each coarse
block (see Figure \ref{schematic_intro} for coarse and fine grid illustration).
For each coarse block $\omega_i$ (or $K$)
and time interval $(T_{n-1},T_n)$,
we identify local multiscale basis functions
$\phi_j^{n,\omega_i}$
($j=1,...,N_{\omega_i}$) and
seek the solution in the span of these basis functions.
For problems with the scale separation, a limited number
of degrees of freedom is sufficient.
For more
complicated heterogeneities as those that appear in many real-world
applications, one needs a systematic approach to find the additional
degrees of freedom.
In each coarse grid (space-time), first,
we construct the snapshot space,
$V_{\text{snap}}^{n,\omega_i}=\text{span}\{\psi_j^{n,\omega_i}$\}.
 The
choice of the snapshot space depends on the global discretization and
the particular application \cite{chung2016adaptive}.
Each snapshot can be constructed, for example, using random boundary
conditions or source terms~\cite{randomized2014}, which allows avoiding
the computations of all snapshot solutions.
Details of snapshot spaces are presented in Appendix \ref{sec:app1}.

Once we construct the snapshot space $V_{\text{snap}}^{n,\omega_i}$,
we compute the offline space,
which is a principal component subspace of the snapshot space.
The offline space
contains important degrees of freedom as first few basis functions.
In this paper, our goal is to develop a Bayesian framework,
which adaptively adds new basis functions to very few
initial basis functions.

The offline space construction is discussed in Appendix \ref{sec:app1}
and in Section \ref{sec:numerics} for some examples.
We denote the offline space by $V_{H,\text{off}}^{n,\omega_i}$
for a generic domain
$\omega_i$ and time interval  $(T_{n-1},T_n)$
with elements of the space denoted $\phi_l^{n,\omega_i}$.
 The offline space is constructed by performing a spectral decomposition
in the snapshot space.
By selecting the dominant eigenvectors
(corresponding to the smallest eigenvalues),
we choose the elements of the offline space \cite{chung2016adaptive}.
 The choice of the spectral problem is important for
the convergence and is derived from the analysis as it is described below.
The convergence rate of the method is proportional to $1/\Lambda_*$.
Here,
$\Lambda_*$ is the smallest eigenvalue among all coarse blocks whose
corresponding
eigenvector is not included in the offline space.
Our goal is to select the local spectral problem
to remove as many small eigenvalues as possible so that we can obtain smaller
dimensional coarse spaces to achieve a higher accuracy.

\section{Bayesian formulation}
\label{sec:bayes}

We seek the solution in each coarse time interval
$(T_{n-1},T_n)$
\begin{equation*}
  u^{n}_H(x,t)=\sum_{i,j}\beta_{i,j}^n \ \phi_i^{n,\omega_j}(x,t) ,
\end{equation*}
where $\beta_{i,j}^n$'s are
defined in each computational time interval and $\phi_i^{n,\omega_j}(x,t)$
are basis functions. We will choose  $\phi_i^{n,\omega_j}$
 to be time-independent, in general.
in this paper.

For further description, we introduce some notations. First,
we select some basis functions in each selected subdomain, which will
be used  as {\it permanent} multiscale basis functions.
These are first dominant modes
(typically a few basis functions) in each coarse domain
and denoted
by
\[
\phi_i^{n,\omega_j}(x,t) \ - \ \text{permanent basis functions}.
\]
Furthermore, we will select basis functions
from the rest of the space, which we denote by
\[
\phi_{i,+}^{n,\omega_j}(x,t)  \ - \ \text{the rest of basis functions}
\]
in order to distinguish from the first few basis functions.
The solution is sought as
\[
 u^{n}_H(x,t)=\sum_{i,j}\beta_{i,j}^n \ \phi_i^{n,\omega_j}(x,t)
+ \sum_{i,j}\beta_{i,j,+}^n \ \phi_{i,+}^{n,\omega_j}(x,t).
\]

{\bf ``Fixed'' multiscale solution with permanent basis functions.}
In our approaches, we will use
``fixed'' multiscale solution computed using
permanent basis function. We denote this solution by
\begin{equation}
  \label{eq:fix}
  u^{n,\text{fix}}_H(x,t)=\sum_{i,j}\gamma_{i,j}^n \ \phi_i^{n,\omega_j}(x,t) ,
\end{equation}
where $\gamma$'s are computed from
\[
A(\gamma^n)= R^{n,\text{fix}}_v(u_{H}^{n}(x,t),u_{H}^{n}(x,T_{n-1}^{+}),u_{H}^{n-1}(x,T_{n-1}^{-}))=0,
\]
where $v$ are offline basis functions.

Next, we describe some general steps of our algorithms,
which are used to define posterior functionals.

\begin{itemize}

\item Problem under consideration
\begin{equation}
\label{eq:main_alg_1}
{\partial u \over \partial t} = \mathcal{L}(u).
\end{equation}

\item We seek the solution in the time interval
$(T_{n-1},T_n)$, as
\[
 u^n_H(x,t)=\sum_{i,j}\beta_{i,j}^n \ \phi_i^{n,\omega_j}(x,t)
+ \sum_{i,j}\beta_{i,j,+}^n \ \phi_{i,+}^{n,\omega_j}(x,t).
\]
Here,  basis functions $\phi_i^{n,\omega_j}(x,t)$  are kept fixed
and $\phi_{i,+}^{n,\omega_j}(x,t)$
are selected based on the indicator $\mathcal{I}$
for basis index and the indicator $\mathcal{J}$ for subdomains.
The indicators are defined later and represent the indices of basis
functions and subdomains that are selected in simulations.

\item Denote the global residual at time $(T_{n-1},T_n)$
by $\mathcal{R}^n$ and the residual for the subdomain
$\omega_j$ by $\mathcal{R}^n_{\omega_j}$. We will consider
several possible residuals:
$ R^{n,\text{fix}}_v(u_{H}^{n,\text{fix}}(x,t),u_{H}^{n,\text{fix}}(x,T_{n-1}^{+}),u_{H}^{n-1,\text{fix}}(x,T_{n-1}^{-}))$;
and
$ R^{n}_v(u_{H}^{n}(x,t),u_{H}^{n}(x,T_{n-1}^{+}),u_{H}^{n-1}(x,T_{n-1}^{-}))$.
In all cases, $v$ is assumed to be in a large dimensional
snapshot space.  $\mathcal{R}^n$
and  $\mathcal{R}^n_{\omega_j}$ are vectors corresponding to these residuals.

\item First, we select $N_{\omega}$ subdomains, where
multiscale basis functions will be added. We denote by
$\alpha^{\omega_k}=\|\mathcal{R}^n_{\omega_k}\|/\|\mathcal{R}^n\|$,
where $\mathcal{R}^n$ is the global residual vector and
$\mathcal{R}^n_{\omega_k}$ is the local residual vector in
$\omega_k$ (as mentioned earlier).
Furthermore, we denote by
\[
\widehat{\alpha}^{\omega_k}={{\alpha}^{\omega_k}\over \sum_j {\alpha}^{\omega_j}}
N_{\omega},
\]
where $N_\omega$ is the desired number of average subregions
that is the user's choice and depends on available computer resources.
With probability $\widehat{\alpha}^{\omega_k}\wedge 1$, we select the region
$\omega_k$, i.e., $\mathcal{J}_k=1$ if the region is selected.

\begin{remark}
  We note that we can compute the residual only at few locations to evaluate the
prior probability distributions.
\end{remark}

\item We use
{\it permanent} basis functions $\phi_i^{n,\omega_j}$ in every region
and select additional
$N_{\text{basis}}$ basis functions (as above) using the residual.
The number of basis,
$N_{\text{basis}}$,
 is the user's choice and depends on available computer resources.
For each $\omega_j$, we compute the correlation coefficient
\[
\text{corrcoeff}(\mathcal{R}^n_{\omega_j}, \phi_{k,+}^{n,\omega_j})=\alpha_{k,+}.
\]
We normalize these $\alpha_{k,+}$ (we keep the same notation
as for subdomain indices) so that on average we have
$N_{\text{basis}}$ basis functions.
\[
\widehat{\alpha}_{k,+}= {\alpha_{k,+}\over \sum \alpha_{i,+}} N_{\text{basis}}.
\]
I.e., in a prior distribution,
we choose $k$th basis with probability $\widehat{\alpha}_{k,+}\wedge 1$.
Since $\widehat{\alpha}_{k,+}$ are used as a prior information,
one can use the residual at a few fixed locations
or approximate residuals to compute this prior.

\item We will consider several posteriors.

 {\bf Posterior around fixed solution.} In this case,
we will use the solution $[0,T]$ computed using permanent
basis to sample realizations. We sample
\begin{equation}
  \label{eq:post1}
  \begin{split}
P(\beta_{+}^{1:N_T},\mathcal{I}^{1:N_T},\mathcal{J}^{1:N_T}|u^{1:N_T,\text{fix}}_H(x,t))
\sim
P(u^{1:N_T,\text{fix}}_H(x,t)|\beta_{+}^{1:N_T}(\mathcal{I}^{1:N_T},\mathcal{J}^{1:N_T}))\\
\pi(\beta_{+}^{1:N_T}|\mathcal{I}^{1:N_T},\mathcal{J}^{1:N_T})
\pi(\mathcal{I}^{1:N_T},\mathcal{J}^{1:N_T}),
\end{split}
  \end{equation}
where superscript $1:N_T$ refers to
the whole time interval of the sampled quantities,
\begin{equation}
\label{eq:post_sigma}
\begin{split}
P(u^{1:N_T,\text{fix}}_H(x,t)|\beta_{+}^{1:N_T}(\mathcal{I}^{1:N_T},\mathcal{J}^{1:N_T}))\sim
 \exp\left( - {\| R^{n,\text{fix}}_v\|^2\over \sigma_L^2} \right)
\end{split}
\end{equation}
and
\[
\pi(\beta_{+}^{1:N_T}|\mathcal{I}^{1:N_T},\mathcal{J}^{1:N_T})=
\exp\left(-{\|\beta_{+}^{1:N_T}\|^2\over\sigma_2^2}\right)
\]
Here, $R^{n,\text{fix}}_v$ is the residual computed using
$u^{1:N_T,\text{fix}}_H(x,t)$, $\|\cdot\|$ is a discrete $l^2$ norm,
and $\sigma_L$ represents the precision of the numerical approximation.
Here, $\sigma_2$ is the variance for the prior and, in our numerical
simulations, we ignore the prior for $\beta_{+}^{1:N_T}$, and so,
one can assume that $\sigma_2$ is a large number.
In fact, the posterior is computed around the residual corresponding to
the fixed solution. This posterior will be used in the example
for parabolic equations (Sections \ref{sec:numerics_parabolic1} and
\ref{sec:numerics_parabolic2}).


{\bf Posterior around previous time.}
In the second case, we will sample at each time interval based
on the previous time sampled solution. In this case, we
re-compute the coefficients corresponding to fixed
basis functions.
\begin{equation}
\label{eq:post2}
\begin{split}
P(\beta^{n+1},\mathcal{I}^{n+1},\mathcal{J}^{n+1}|u^{n}_H(x,t)))
\sim P(u^{n}_H(x,t))|\beta^{n+1}(\mathcal{I}^{n+1},\mathcal{J}^{n+1}))\\
\pi(\beta^{n+1}|\mathcal{I}^{n+1},\mathcal{J}^{n+1})
\pi(\mathcal{I}^{n+1},\mathcal{J}^{n+1})
\end{split}
\end{equation}
\begin{equation}
  \begin{split}
    P(u^n_H(x,t))|\beta^{n+1}(\mathcal{I}^{n+1},\mathcal{J}^{n+1}))\sim
    \exp\left(- {\| R^{n+1}_v\|^2\over\sigma_L^2}  \right)
\end{split}
  \end{equation}
Here, $R^{n+1}_v$ is the residual computed using
$u^{n}_H(x,t)$ as an initial condition and permanent basis functions
in the current time interval $(T_n,T_{n+1})$,
$\|\cdot\|$ is a discrete $l^2$ norm,
and, again, $\sigma_L$ represents the precision of the numerical approximation.
In fact, the posterior is computed around the residual corresponding to
the solution at previous time step. This posterior will be used in the example
for wave equations (Sections \ref{sec:numerics_wave}).

{\bf Posterior using fixed and previous time solutions.}
One can condition the multiscale solution at the fixed solution
$u^{1:N_T,\text{fix}}_H(x,t)$ and $u^{n}_H(x,t)$. In this case,
one can seek a posterior distribution for
\[
P( u^{n+1}_H(x,t)| u^{n}_H(x,t), u^{1:N_T,\text{fix}}_H(x,t)).
\]
Other more general posteriors can also be setup.

\end{itemize}


\subsection{Sampling from the posterior distribution}

We will consider two sampling algorithms. In the first sampling
algorithm (we call ``sequential sampling''), we will use
samples from the prior distribution and use them to obtain
samples of $\beta$'s. In this process, in the first step,
we sample $\mathcal{I}$ and $\mathcal{J}$ from the prior
based on the residual
and use them to compute the samples of $\beta$. In our numerical
experiments, we will show the results for
$E(\beta^{n+1}|\mathcal{I}^{n+1},\mathcal{J}^{n+1})$.

In the second algorithm, we will perform full posterior sampling
(we call ``full sampling''), where we will sample both the indices
($\mathcal{I}$, $\mathcal{J}$)
and the coefficients ($\beta$). In this case, we will use Gibbs
sampling \cite{robert2004monte,gilks1995markov}, though one can
also design efficient sampling algorithms
based on Markov chain Monte Carlo methods \cite{robert2004monte,gilks1995markov}.
In Gibbs sampling,
we will compute the probabilities
$\hat{\pi}_i$ for each additional
basis function as
\[
\frac{\hat{\pi}^{n+1}_{i+}}{1-\hat{\pi}^{n+1}_{i+}}=
\frac{\hat{\alpha}^{n+1}_{i+}}{1-\hat{\alpha}^{n+1}_{i+}}*\mathcal{F},
\]
where $\mathcal{F}$ is a multiplication factor that depends how much
a change in one basis function
(i.e., adding a basis function)
will affect the residual change.
This factor also depends on the model size. To penalize
adding linearly dependent basis functions, we modify the prior
distribution $\pi(\mathcal{I}, \mathcal{J})$
with a multiplicative constant, which takes into account the linear
dependency factor.
This multiplicative factor is a product of the singular values of the
matrix consisting of an inner product of basis functions such
that the resulting $\sigma^2\log(\mathcal{F})$ reduces to the difference
 between
two residuals.
Note that in this step, it is important to have a computable residual.

\subsection{Dynamic data}

In our proposed approach, the threshold $\sigma_L$ (e.g.,
in (\ref{eq:post_sigma}))
represents the accuracy of our approximation. If we let
$\sigma_L$ approach to zero, our sampling will increase the
number of basis functions and the solution will converge to the
fine-grid solution. However, the main idea of our approach
is to allow less accurate solutions.
The
 threshold in (\ref{eq:post_sigma}) can depend on the additional
dynamic data.
We can add the time-dependent data by including additional
terms in the posterior distribution. For example, if
we obtain the measurements of the solution at some locations,
this can be added as an additional multiplicative term
in the posterior (e.g., in (\ref{eq:post1})) in the form of
\[
\exp\left(-{ \|\mathcal{D}(u_H^{n+1})-\mathcal{D}_{obs}\|^2\over \sigma_d^2}  \right),
\]
where $\mathcal{D}(u_H^{n+1})$ is the observation that depends
on the solution and
$\mathcal{D}_{obs}$ is the associated observed data.
We note that $\sigma_d$ is a factor in choosing the accuracy of
our approximation. For example, if the measurement accuracy $\sigma_d$
is large,
one can choose $\sigma_L$ also to be large.

\section{Numerical examples }
\label{sec:numerics}

In our numerical results, we will compare two sampling
approaches. In the first approach,
the basis functions will be selected from the prior distribution
and we will show the mean of the conditional distributions.
 In the second
approach, we will apply Gibbs sampling to sample the realizations.
Our numerical results will show that both approaches provide a
good accuracy
(in a statistical sense) and the Gibbs sampling is more accurate compared
when sampling from the prior distribution. Moreover,
we observe a fast convergence
when using Gibbs sampling.

We will consider several problems. In the first example, we consider a
parabolic problem and two different discretizations. In the second example,
we will consider the wave equations. In each example, we will specify
the residual that is used in the sampling. We note that
it is important that
this residual provides a good stability for the solution.


\subsection{Mixed formulation for parabolic equations}
\label{sec:numerics_parabolic1}


For the mixed GMsFEM,
the support of multiscale basis functions are
$\omega_E$, which are the two coarse elements sharing a common edge $E$ (see \cite{chung2015mixed} for
details).
In particular, we denote
$\mathcal{E}^{H}$ as the set of all coarse grid edges
and let $N_{e}$ be the total number of coarse grid edges. The coarse grid neighborhood
$\omega_{E}$ of a face $E \in \mathcal{E}^{H}$ is defined as
$$\omega_{E} = \bigcup \{K \in \mathcal{T}^{H} : E \in \partial K \}, \quad i=1,2,\cdots, N_e,$$
which is a union of two coarse grid blocks if $E_i$ is an interior edge (face)
For a coarse edge $E_i$, we write $\omega_{E_i} = \omega_i$ to unify the notations.
Multiscale basis functions are constructed
using local snapshot spaces and local spectral problems
(see \cite{chung2015mixed} for details). Here, we will discuss the residual
which is used in our Bayesian Multiscale Method.

We denote $Q_{H,\text{off}}$ as the space of functions which are piecewise constant on each coarse block. We will use this space to approximate the pressure $u$. For approximating the velocity $\mathcal{V}=-\kappa \nabla u$, we will
construct a multiscale space $V_{H,\text{off}}$ for the velocity by
following the general framework of the GMsFEM.
Next, we use the spaces $Q_{H,\text{off}},V_{H,\text{off}}$ to solve
the problem: to find $u^{n+1}_{H} \in Q_{H,\text{off}}$ and
$\mathcal{V}^{n+1}_{H} \in V_{H,\text{off}}$ such that
\begin{equation}\label{MsEquation}
\begin{split}
\int_{\Omega} \kappa^{-1} \mathcal{V}^{n+1}_{H} \cdot w - \int_{\Omega} \mbox{div}(w) u^{n+1}_{H}= 0, \ \forall w \in V_{H,\text{off}}, \\
\int_{\Omega} (\cfrac{u^{n+1}_{H}-u^{n}_{H}}{\Delta t} + \mbox{div}(\mathcal{V}^{n+1}_{H})) q = \int_{\Omega} f^{n+1}q, \ \forall q \in Q_{H,\text{off}}.
\end{split}
\end{equation}
The residual is defined as
\[
R_w^{n}(\mathcal{V}^{n+1}_{H,+},\mathcal{V}^{n+1,\text{fix}}_H,u^{n+1,\text{fix}}_H)= \int_{\Omega} \kappa^{-1} (\mathcal{V}^{n+1}_{H,+}+\mathcal{V}^{n+1,\text{fix}}_H) \cdot w - \int_{\Omega} \mbox{div}(w) u^{n+1,\text{fix}}_H.
\]
Denote $\mathcal{R}^n$ to be corresponding coordinates
and use discrete $l^2$ norm as the residual norm.

\begin{figure}[H]
\begin{centering}
\includegraphics[scale=0.5]{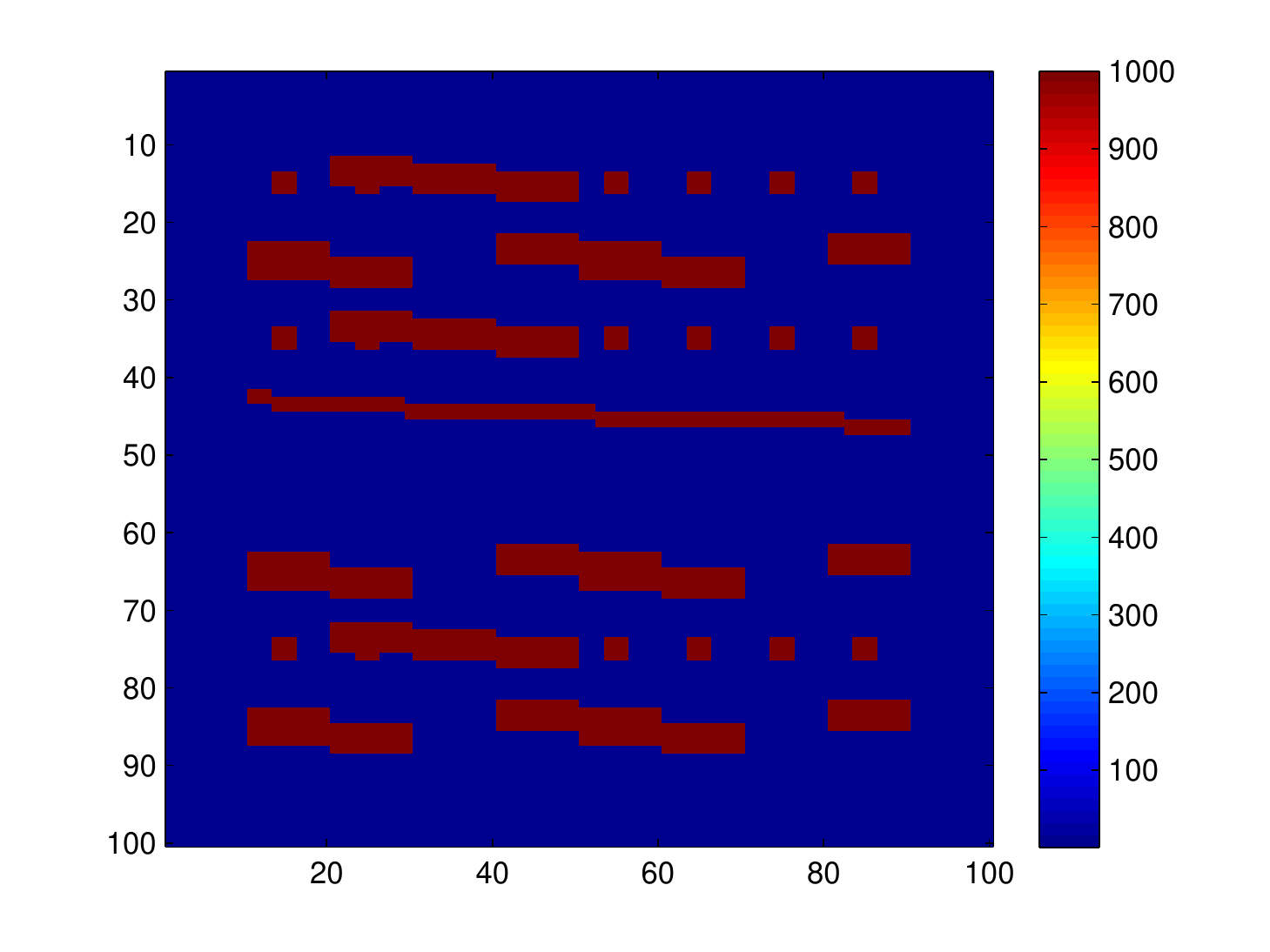}
\par\end{centering}

\protect\caption{The permeability field $\kappa$.}
\label{fig:medium_mixed}
\end{figure}

Next, we will present an example. We use the permeability
field $\kappa$ shown in Figure \ref{fig:medium_mixed} and the contrast of the permeability field $\cfrac{\max{\kappa}}{\min{\kappa}}$ is increasing as  $\cfrac{\max{\kappa}}{\min{\kappa}} = 1000e^{250t}$.
In this example, we find the distribution of the solution
at two different time instants $T=0.01$ and $T=0.02$.
The fine grid is $100\times 100$ and the coarse grid is $10\times 10$.
We use only one permanent basis function per edge to
compute ``fixed'' solution and use our Bayesian framework
to seek additional basis functions by solving small global
problems. In our approach, first using the global residual we define
local regions, where multiscale basis functions are added.
We choose $\sigma_L=1e-3$.
In our example, we fixed these local regions in each time step and choose
$25$\% and $55$\% of the total coarse edges for the first and second time step, where the residual is larger
than a certain threshold. In these coarse blocks,
we apply both sequential sampling and full sampling algorithms.

\begin{figure}[H]
\begin{centering}
\includegraphics[width=1.5in,height=1.5in]{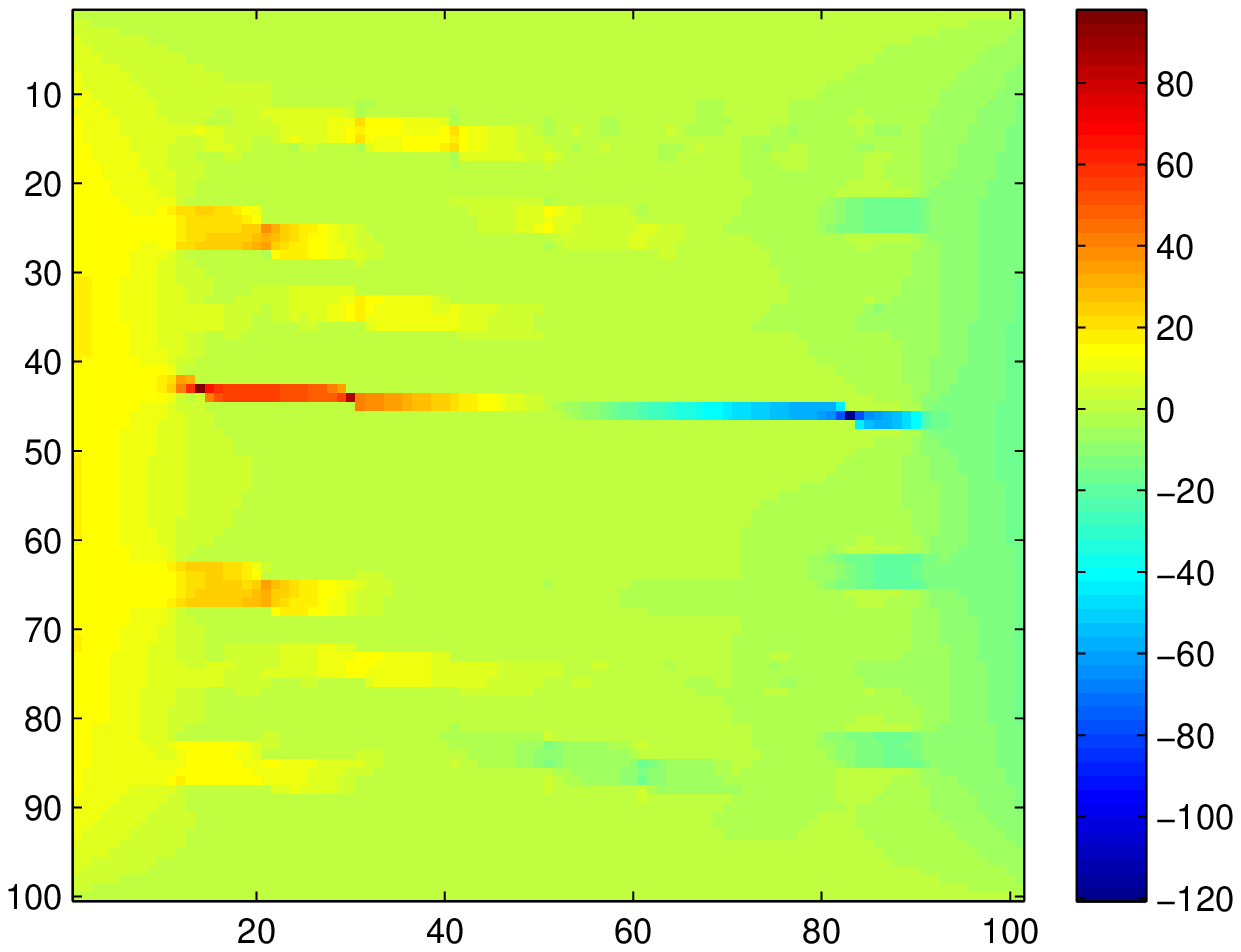}
\includegraphics[width=1.5in,height=1.5in]{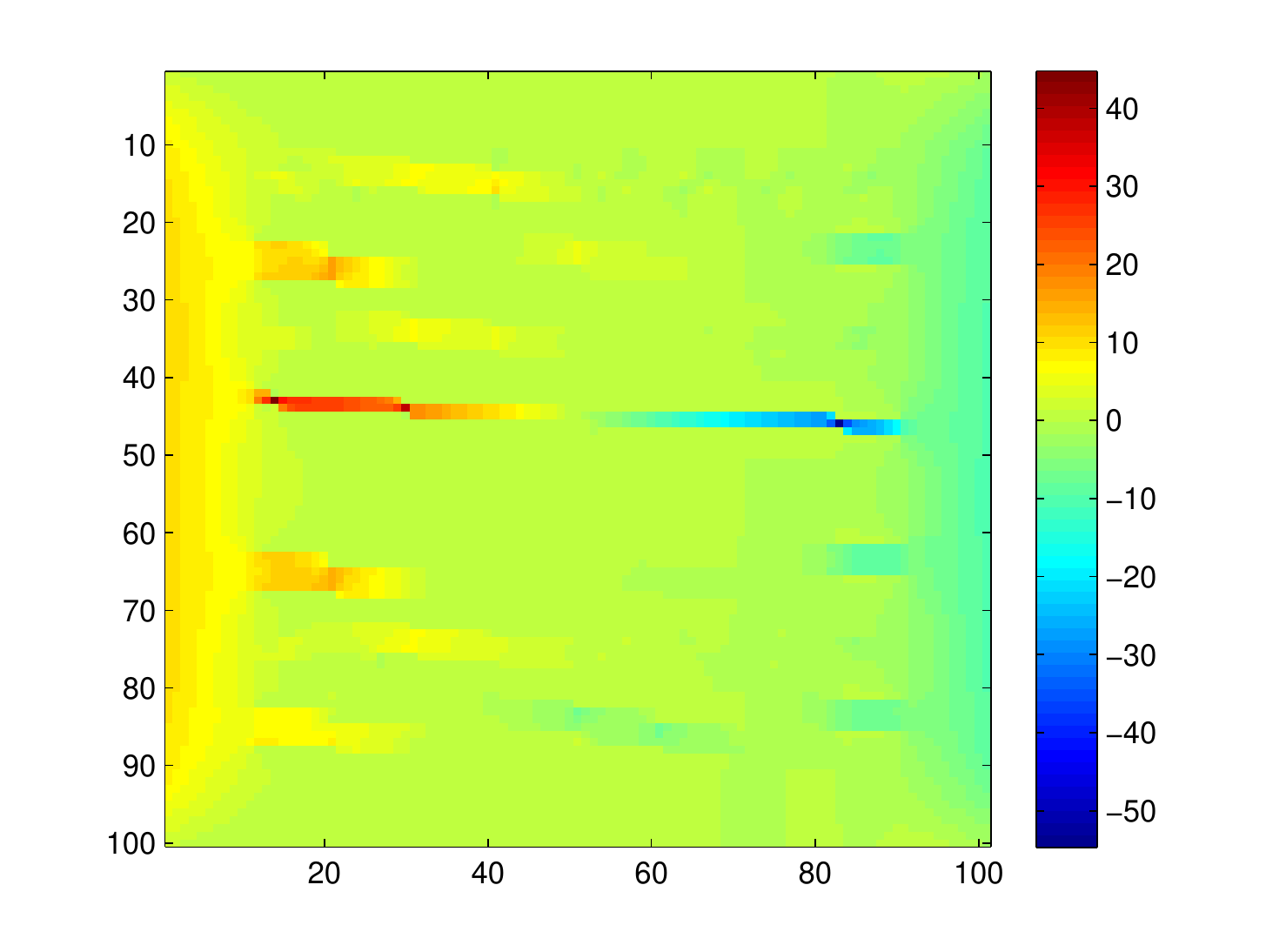}
\includegraphics[width=1.5in,height=1.5in]{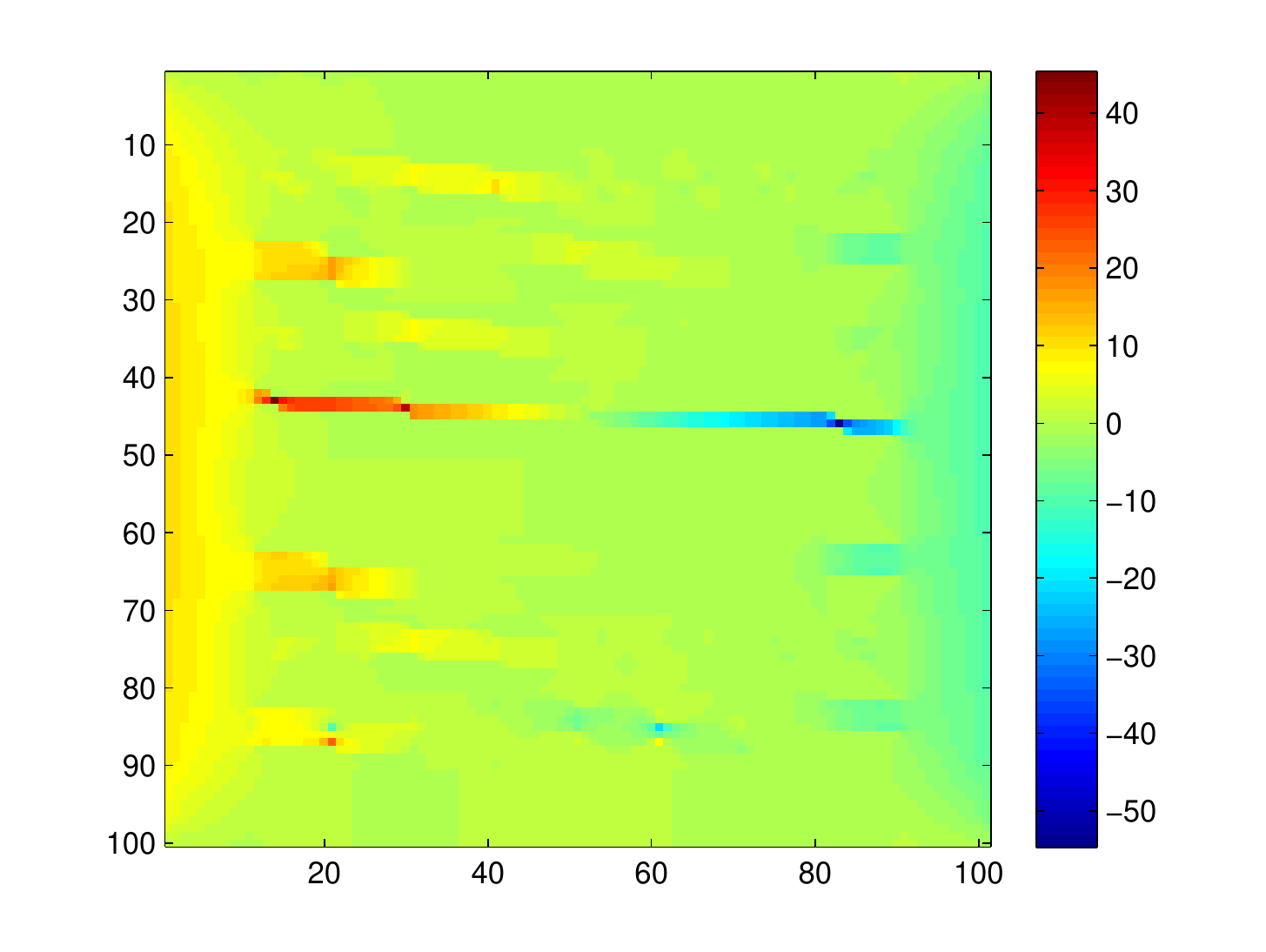}
\par\end{centering}
\protect\caption{Plots of $x_1$-component of the numerical velocity
$\mathcal{V}$ at $T = 0.02$: reference solution (left), mean of sequential sampling (middle), mean of full sampling (right).}
\label{fig:mean_alg1}
\end{figure}

\begin{figure}[H]
\begin{centering}
\includegraphics[scale=0.4]{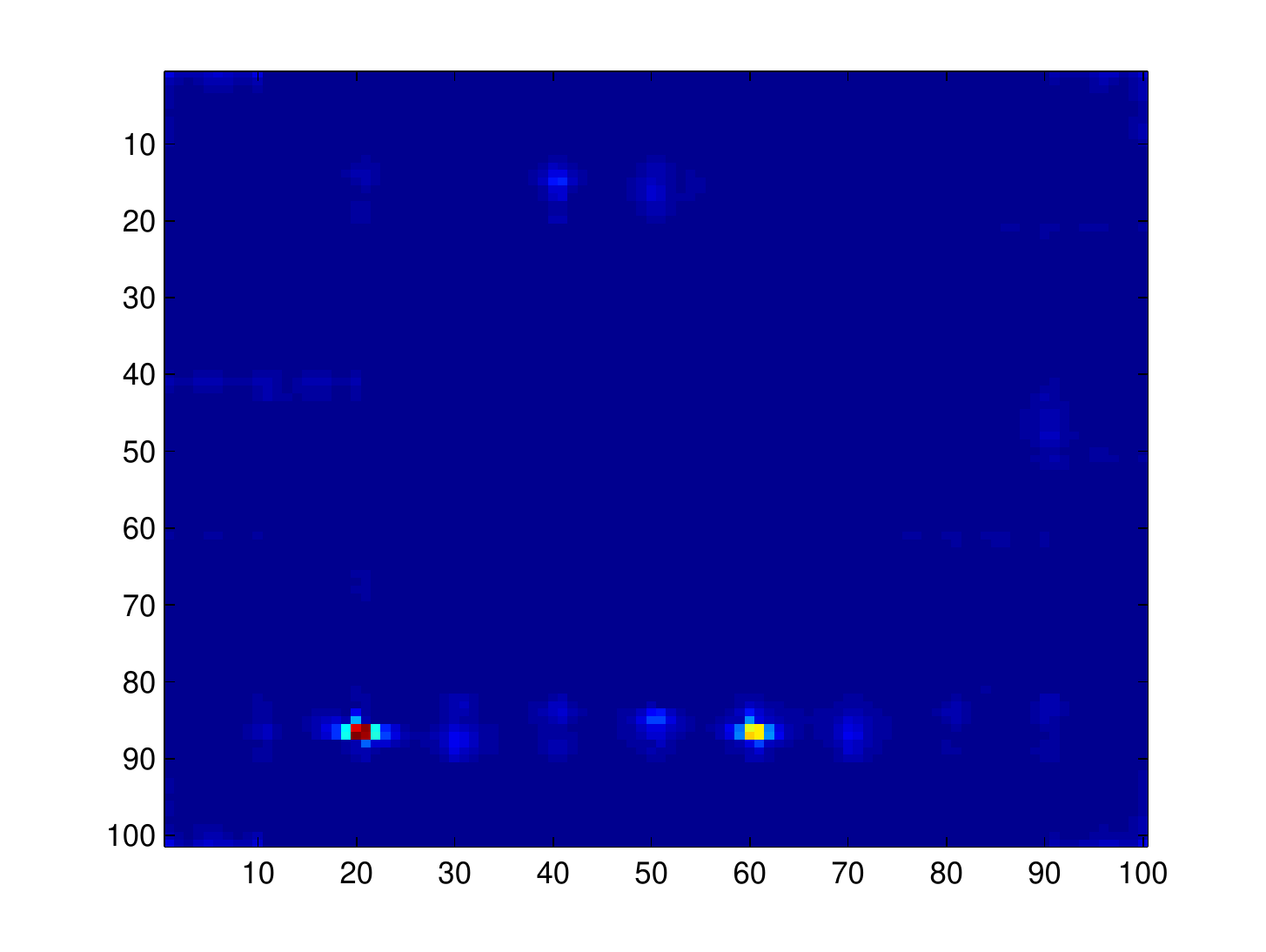}\includegraphics[scale=0.4]{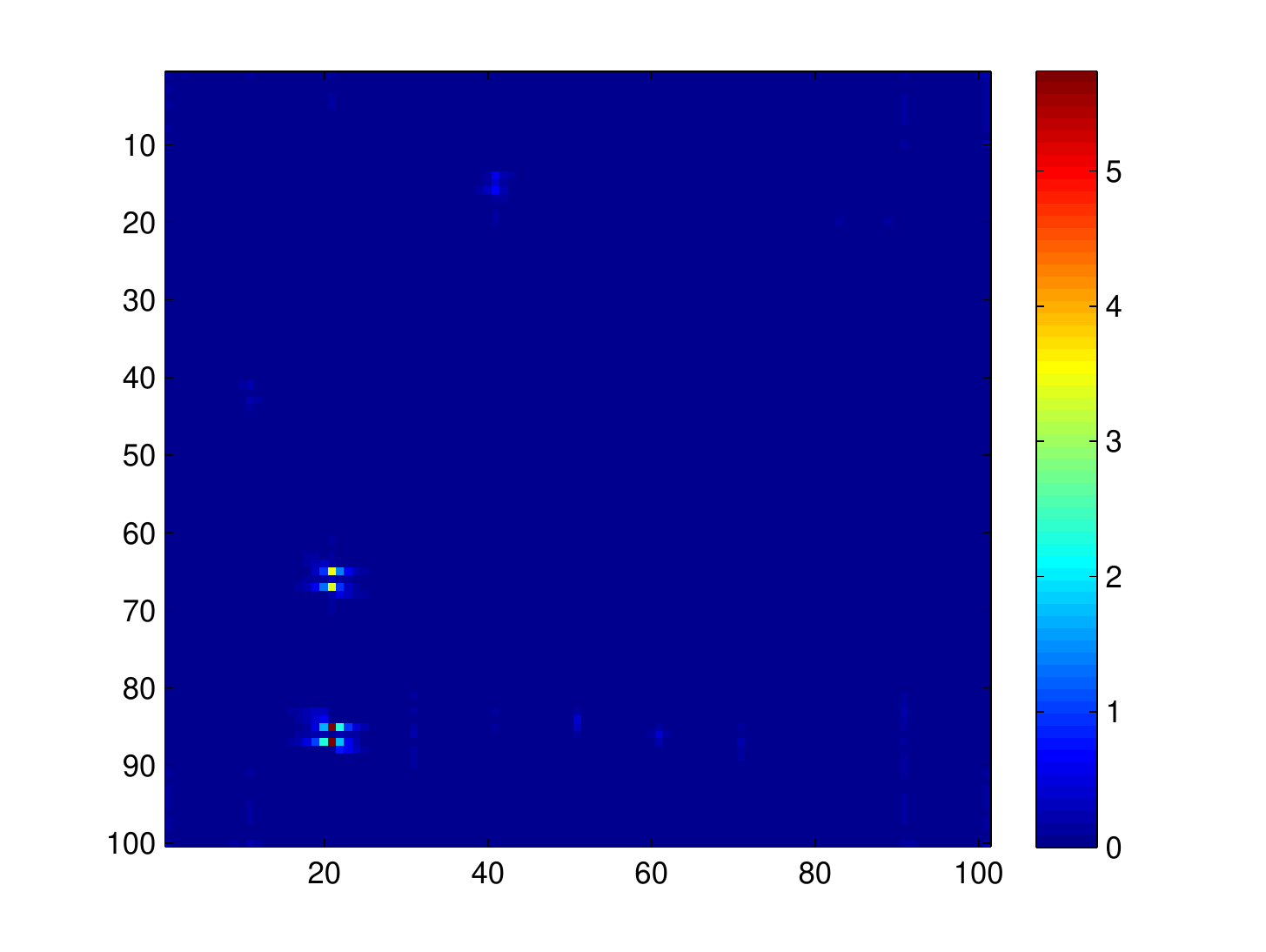}
\par\end{centering}
\protect\caption{Plots of sample standard deviation of $x_1$-component of the velocity at $T = 0.02$: sequential sampling (left), full sampling (right).}
\label{fig:std_alg1}
\end{figure}

\begin{figure}[H]
\begin{centering}
\includegraphics[scale=0.4]{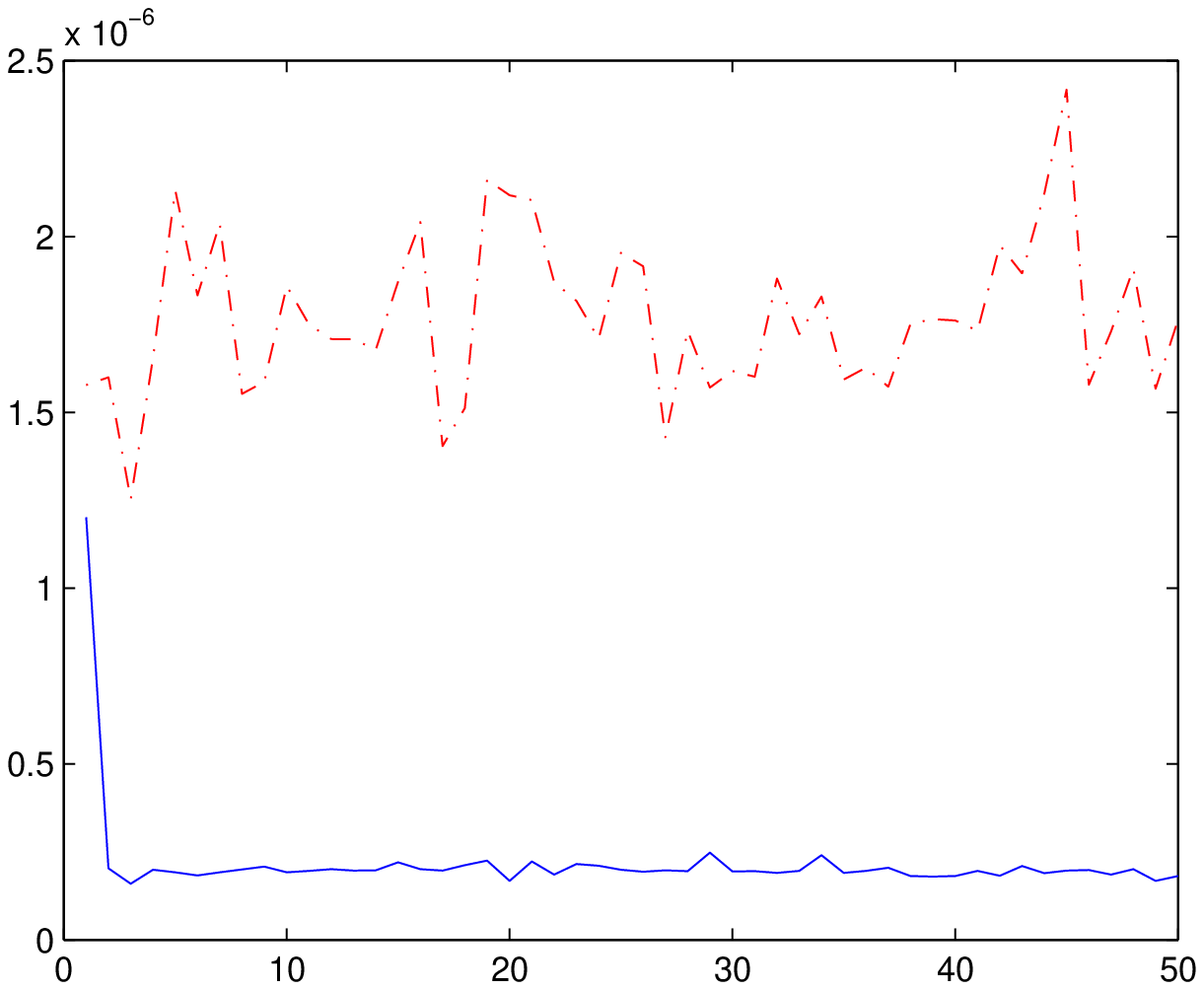}\includegraphics[scale=0.4]{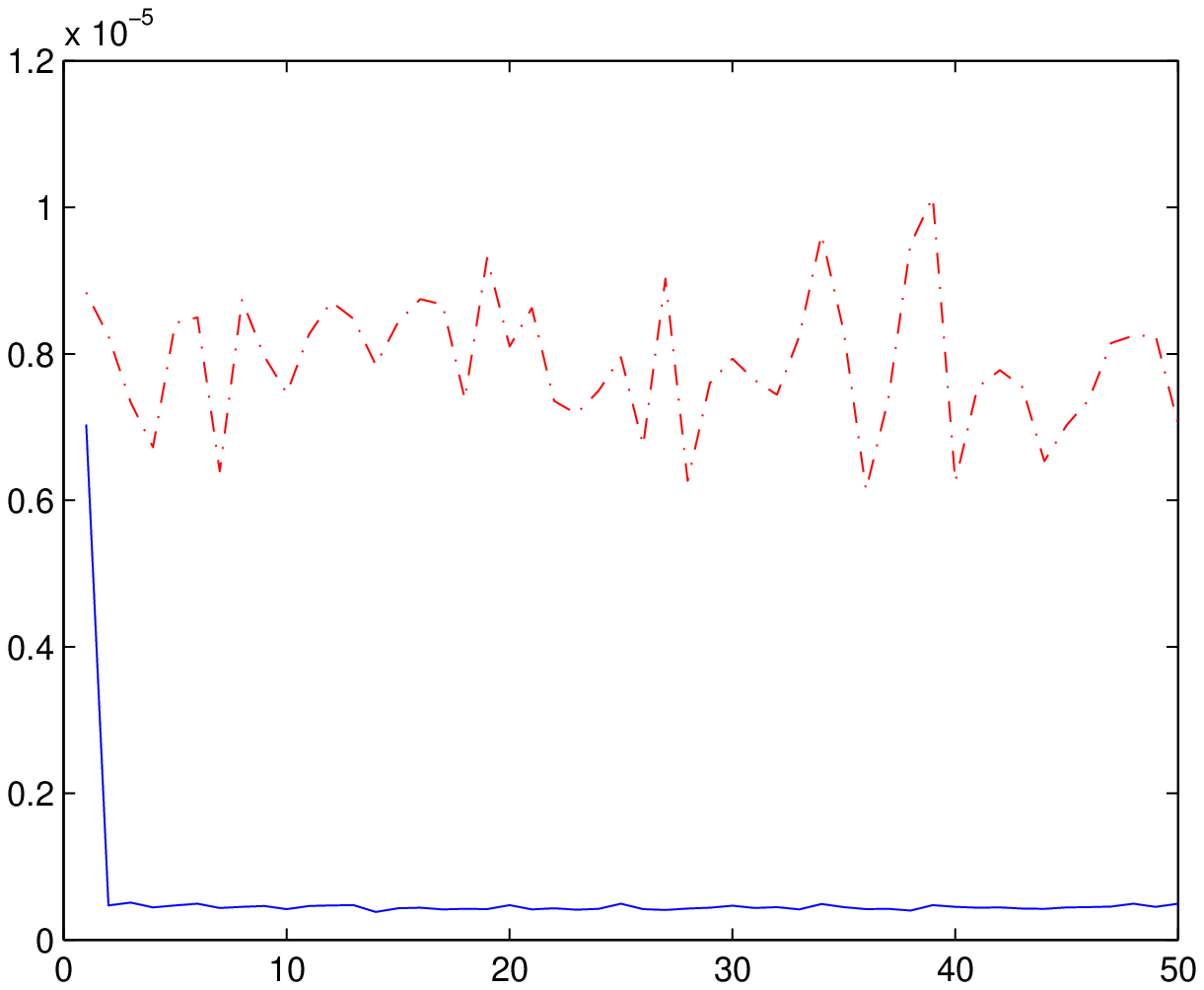}
\par\end{centering}
\protect\caption{Residual vs. samples for
 sequential sampling (red dotted line) and full sampling (blue solid line): at time $T = 0.01$ (left), at time $T = 0.02$ (right).}
\label{fig:residual_Ex1}
\end{figure}

\begin{figure}[H]
\begin{centering}
\includegraphics[scale=0.4]{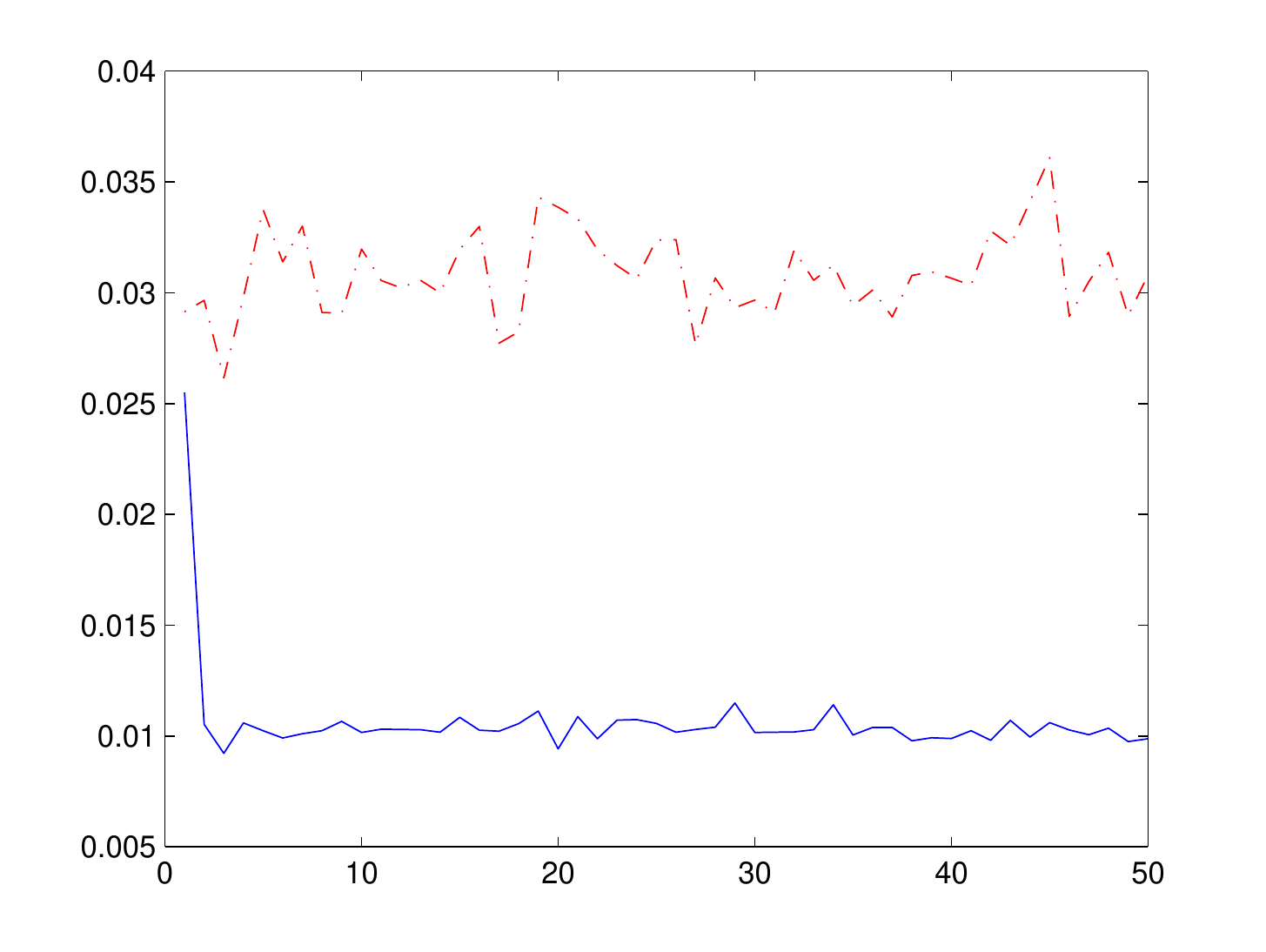}\includegraphics[scale=0.4]{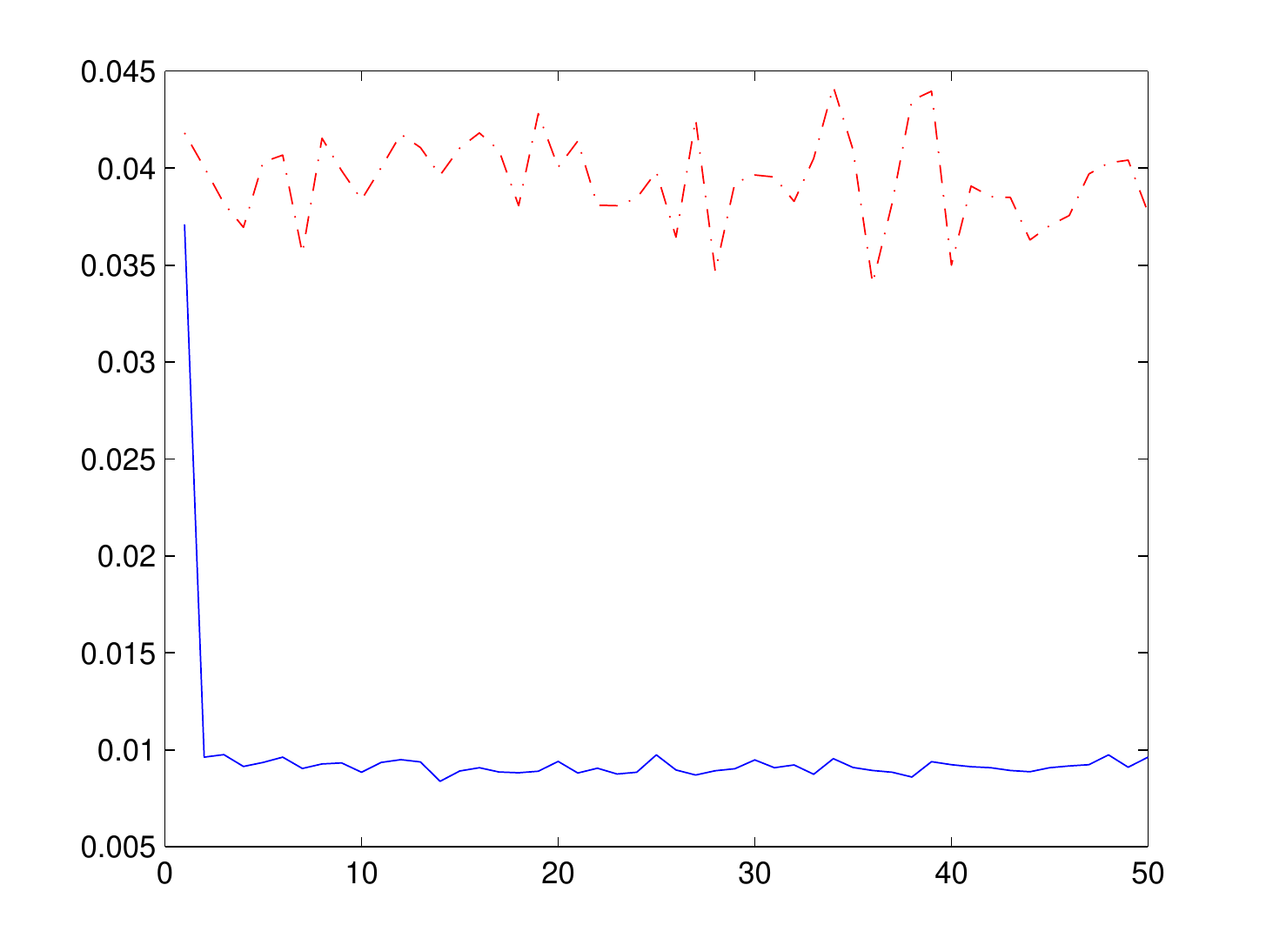}
\par\end{centering}

\protect\caption{$L^2$ vs. sample using sequential sampling (red dotted line) and full sampling (blue solid line): at time $T = 0.01$ (left), at time $T = 0.02$ (right).}
\label{fig:error_Ex1}
\end{figure}

In Figure \ref{fig:mean_alg1}, we depict the mean
solution using the sequential sampling algorithm and full sampling
algorithms. The errors for the mean at $T = 0.02$ are $3.02\%$ for the sequential sampling
and $0.79\%$ for full sampling. We note that full sampling provides a better
result. We observe similar results for $x_2$ component of the solution.
We depict the standard deviation of the solution at each pixel
in Figure \ref{fig:std_alg1}. We observe that the
true solution falls within the limits of the
 mean and the standard
deviations. Next, we show the results across several samples.
For the sequential sampling, we use $20$ realizations and show both
residuals and the errors in Figure \ref{fig:residual_Ex1} and \ref{fig:error_Ex1}. The error is computed as a difference between the solution
and the snapshot solution using the snapshot vectors in the elements,
where the basis functions are updated.
From these figures, we observe that the residuals and errors are
smaller for full sampling compared to sequential sampling. Moreover,
Gibbs sampling stabilizes in a few iterations,
which shows a fast numerical convergence.

\begin{figure}[H]
\begin{centering}
\includegraphics[scale=0.175]{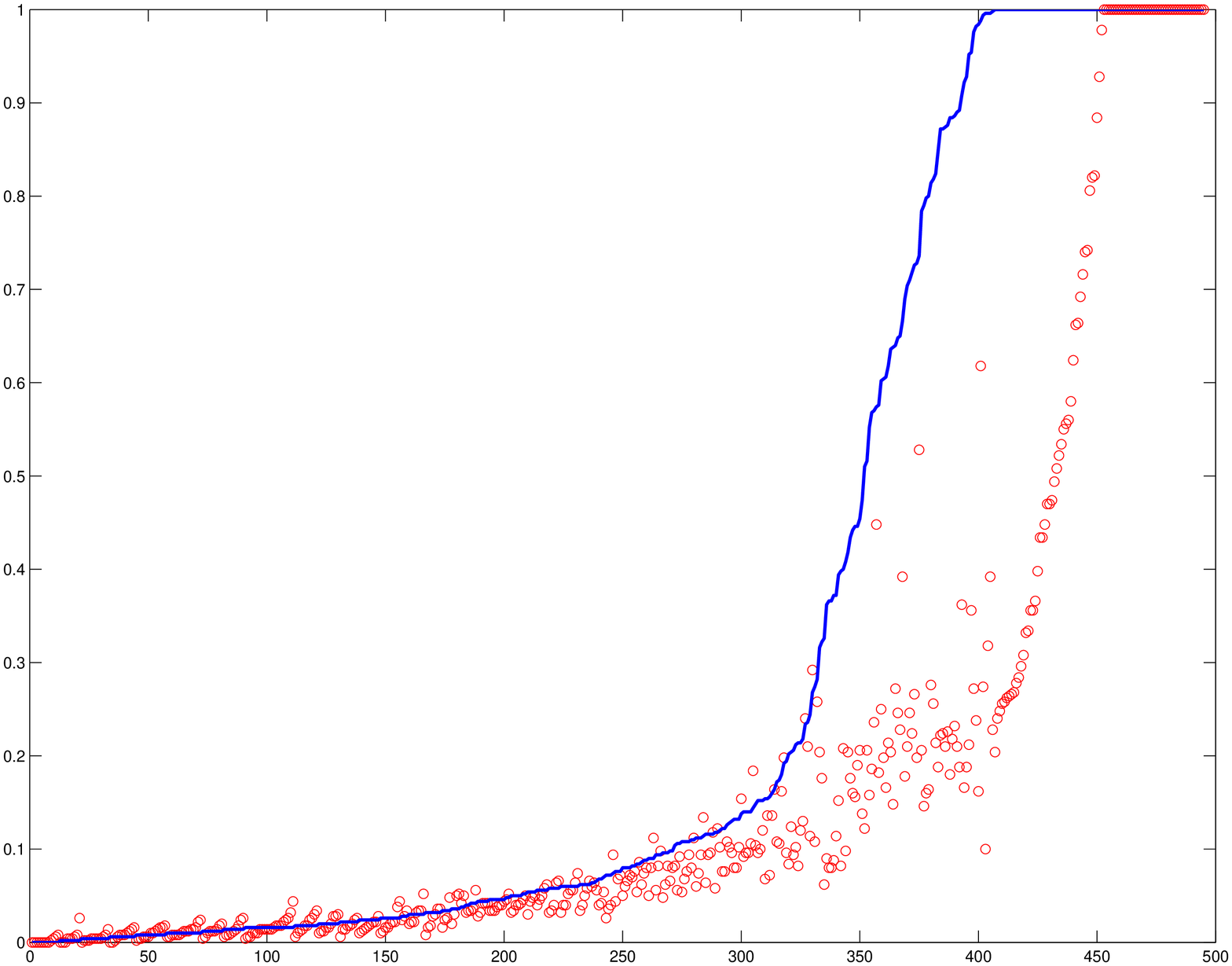}\includegraphics[scale=0.175]{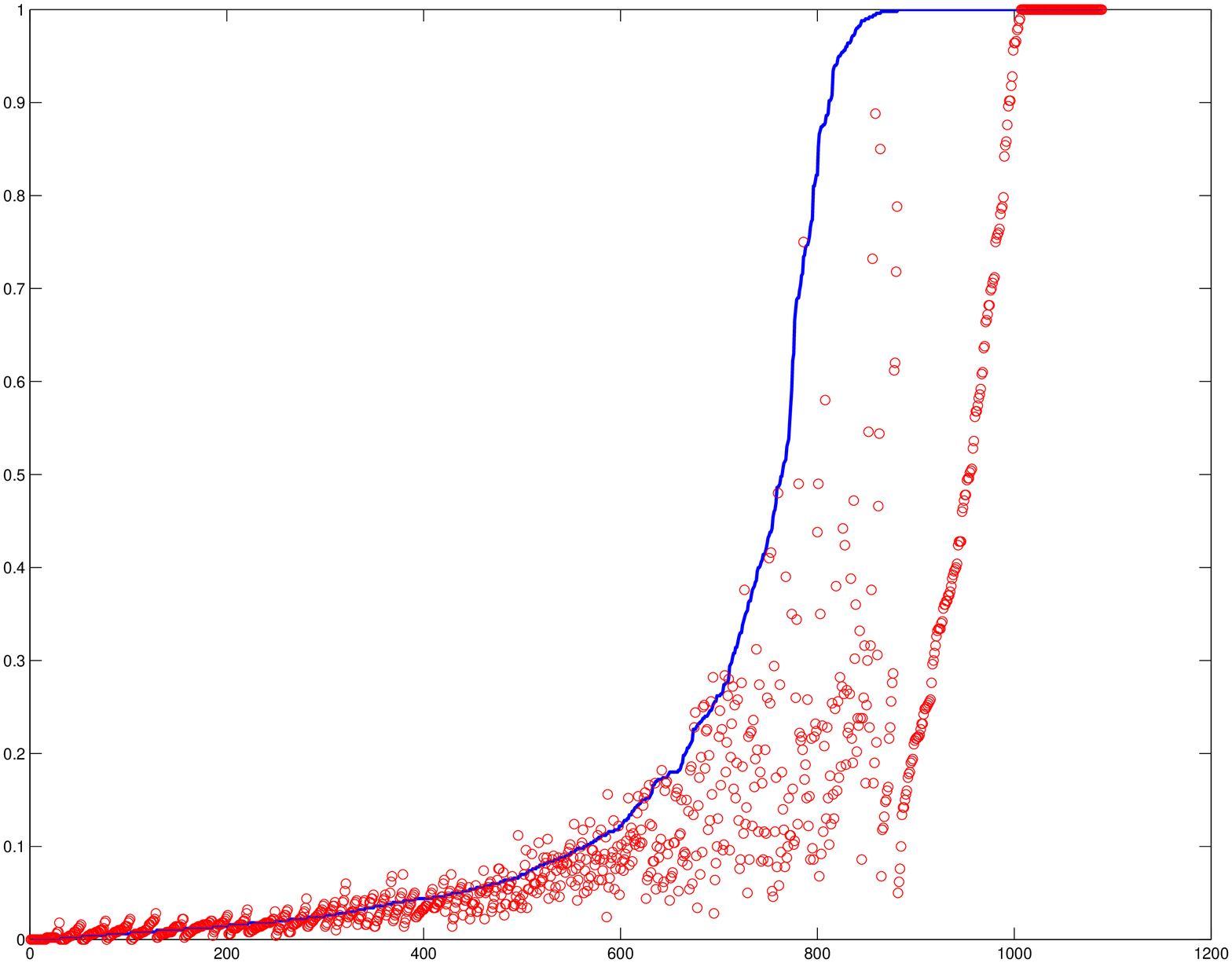}
\par\end{centering}
\caption{History of occurrence probability against basis functions using sequential sampling (red dotted line) and full sampling (blue solid line): at time $T = 0.01$ (left), at time $T = 0.02$ (right).}
\label{fig:mixed_basis_vs_sample}
\end{figure}

\begin{figure}[H]
\begin{centering}
\includegraphics[scale=0.4]{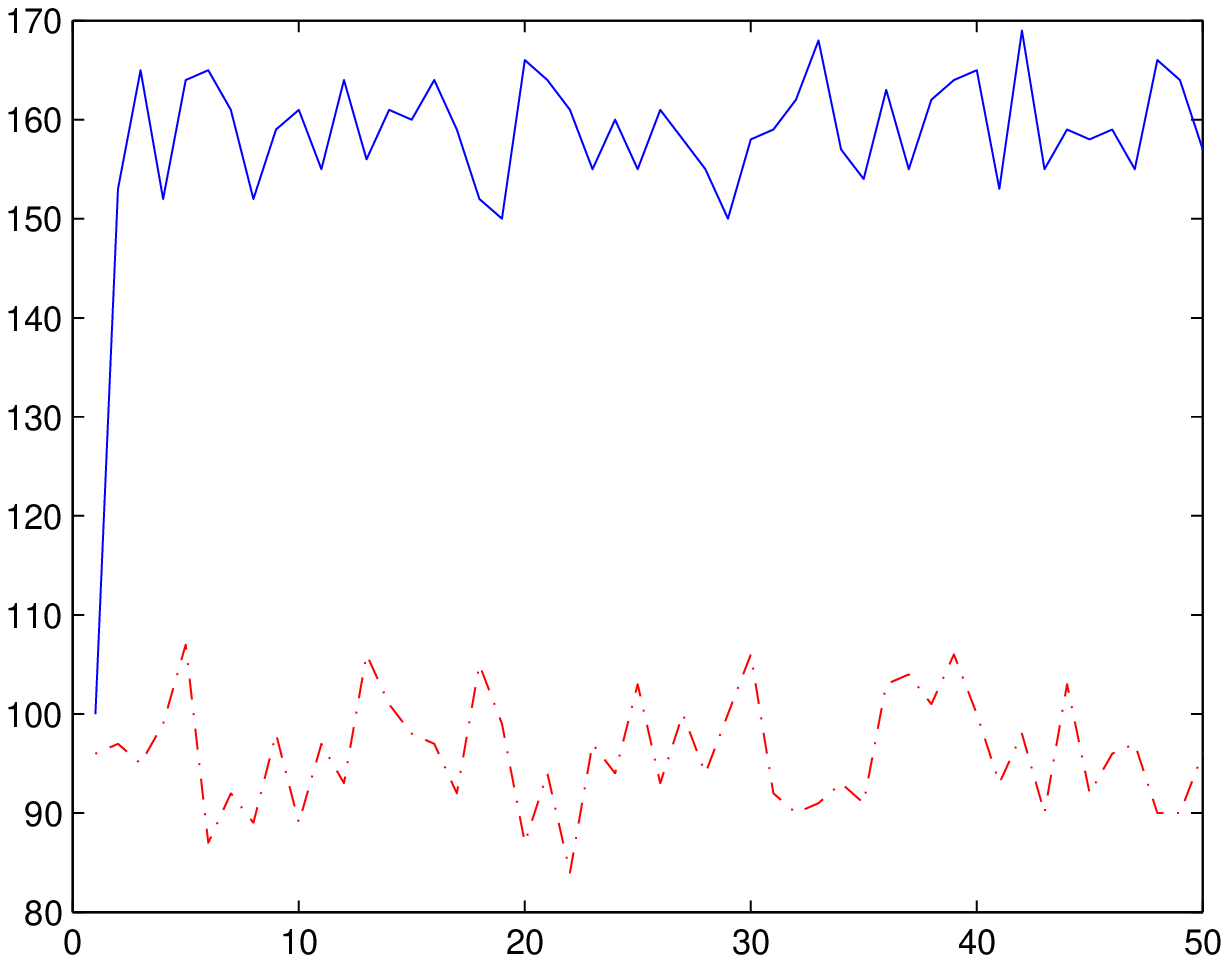}\includegraphics[scale=0.4]{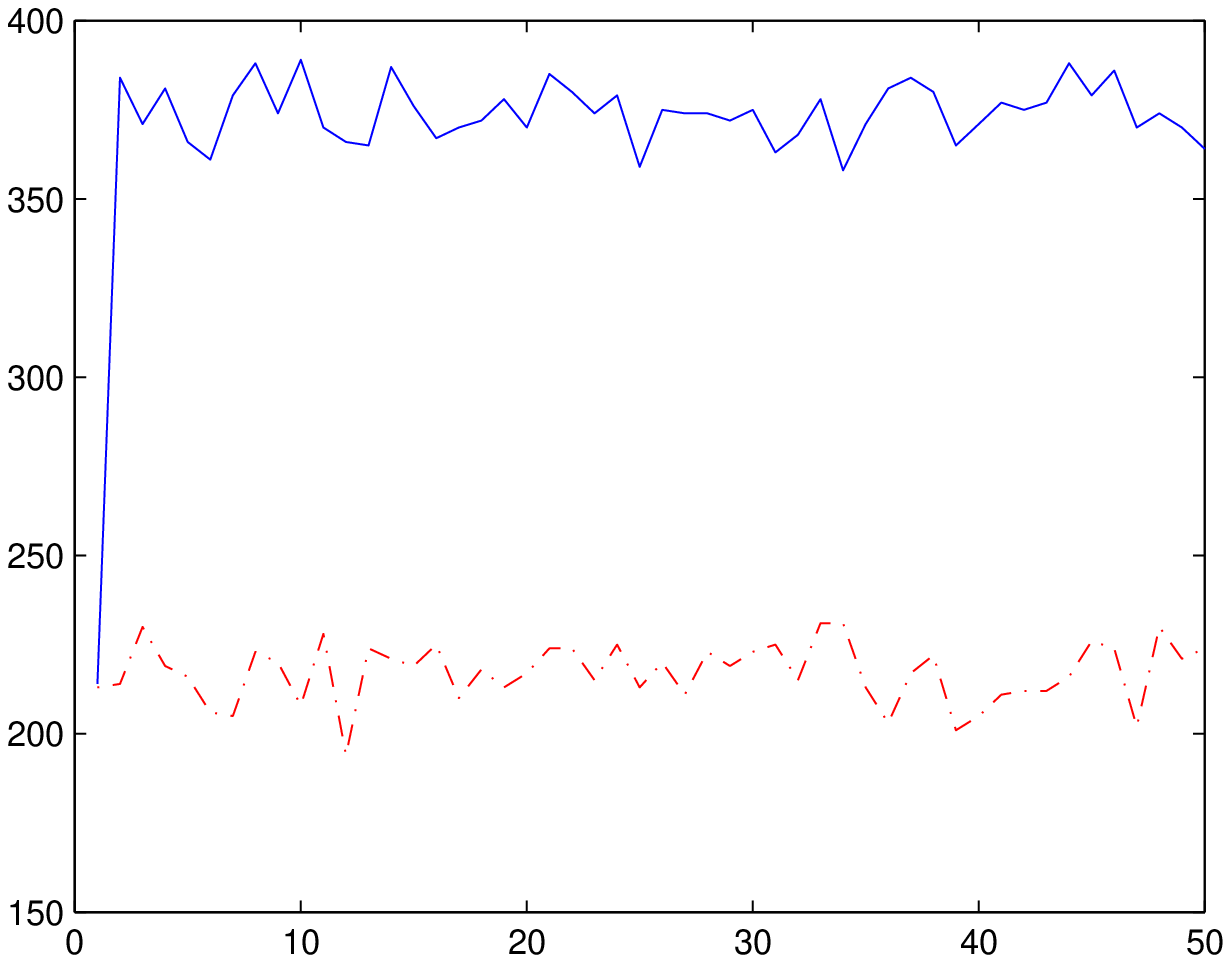}
\par\end{centering}
\protect\caption{History of number of basis functions against sampling process using sequential sampling (red dotted line) and full sampling (blue solid line): at time $T = 0.01$ (left), at time $T = 0.02$ (right).}
\label{fig:mixed_total_basis}
\end{figure}

In Figure \ref{fig:mixed_basis_vs_sample}, we depict
the frequency of the basis functions in both sequential and full sampling.
We have enumarated the basis functions such that the frequency in
full sampling is increasing. We observe from this figure
that the sequential sampling has a similar trend and
we have found the correlation in the frequencies
 between full sampling and partial sampling to be $0.7$.
The shift in sequential sampling is due to the fact that fewer
basis are used in this sampling. We plot the total
number of basis functions
in Figure \ref{fig:mixed_total_basis}. We observe from
this figure that the full sampling requires more
basis functions compared to the sequential sampling.
Moreover, the number of basis functions in full sampling
stabilizes around a certain value, which depends on
$\sigma_L$ (the precision of the residual in (\ref{eq:post_sigma})).

\subsection{Continuous Galerkin formulation for parabolic equations}
\label{sec:numerics_parabolic2}

In this section, we present the Bayesian approach
for the continuous Galerkin formulation.
Multiscale basis functions are obtained from eigenfunctions in the local snapshot space with small eigenvalues in an appropriate local spectral eigenvalue problem (see \cite{chung2016adaptive} for details). We denote the space of multiscale basis functions by
$V_{H,\text{off}}$, in which we seek numerical approximations for the problem: find $u_H^{n+1} \in V_{H,\text{off}}$ such that
\begin{equation}
\int_{\Omega} \cfrac{u^{n+1}_{H}-u^{n}_{H}}{\Delta t} \,  v + \int_{\Omega} \kappa \nabla u^{n+1}_{H} \cdot \nabla v  = \int_{\Omega} f^{n+1}v, \ \forall v \in V_{H,\text{off}}.
\end{equation}
The residual is defined as
\begin{equation}
\begin{split}
R^n_v(u_{+}^{n+1},u_{\text{fixed}}^{n+1},u_{\text{fixed}}^{n}) & = \int_{\Omega} f^{n+1}v - \int_{\Omega} \cfrac{u^{n+1}_{+}+u_{\text{fixed}}^{n+1}-u_{\text{fixed}}^{n}}{\Delta t} \,  v \\
& \quad \quad + \int_{\Omega} \kappa \nabla (u^{n+1}_{+}+u_{\text{fixed}}^{n+1}) \cdot \nabla v.
\end{split}
\end{equation}

We use the permeability
field $\kappa$ as in Figure \ref{fig:medium_mixed}.
We will compare the solutions
at two different time instants $T=0.01$ and $T=0.02$.
The fine grid is $100\times 100$ and the coarse grid is $10\times 10$.
We use $2$ permanent basis functions per coarse neighborhood to
compute ``fixed'' solution and use our Bayesian framework
to seek additional basis functions by solving small global
problems.
In this example, we select $30\%$ of the total
local regions at which residual is the largest
and multiscale basis functions are added.
In these coarse blocks, we apply both sequential sampling
and full sampling.

\begin{figure}[ht!]
\centering
\includegraphics[width=1.5in, height=1.5in]{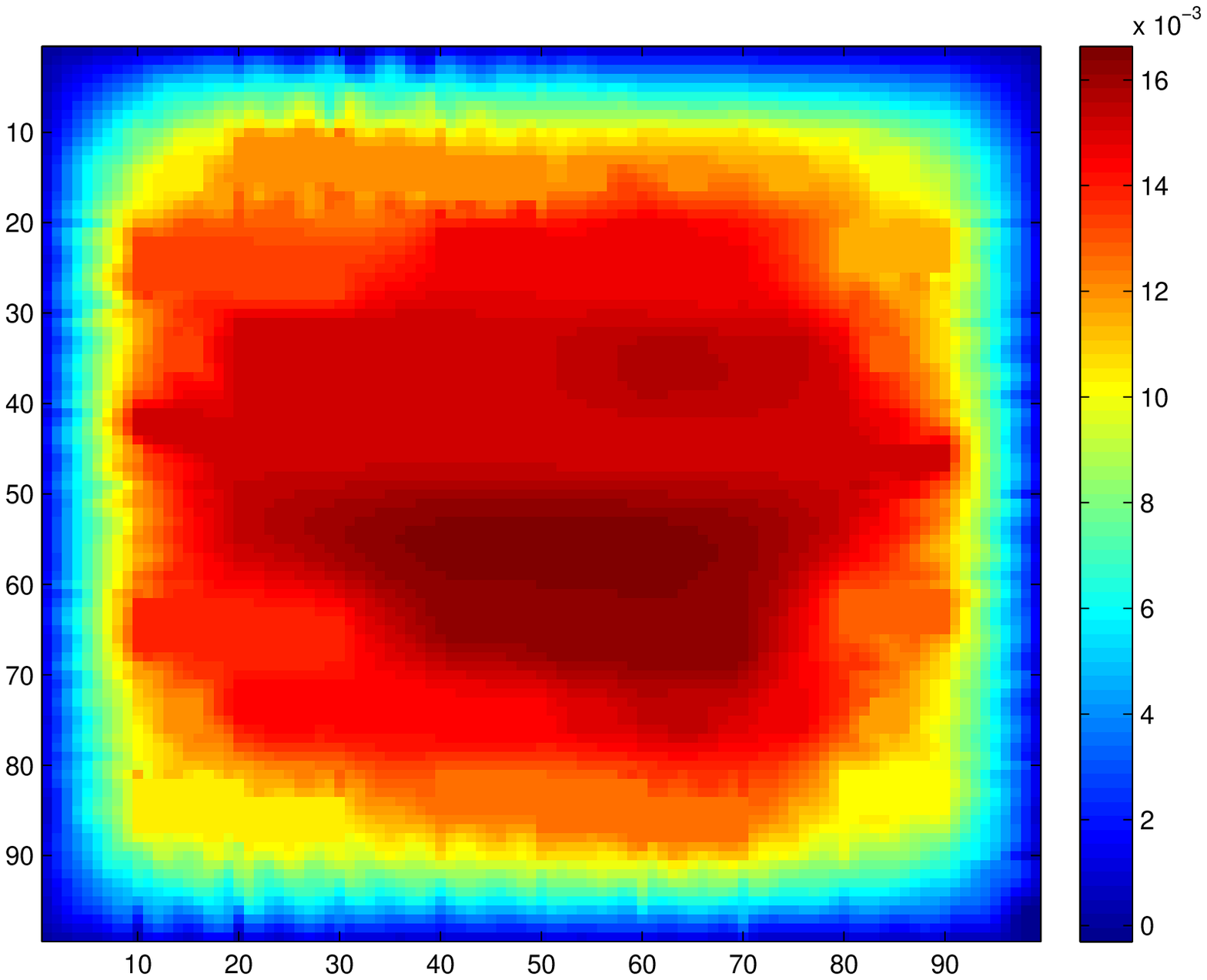}
\includegraphics[width=1.5in, height=1.5in]{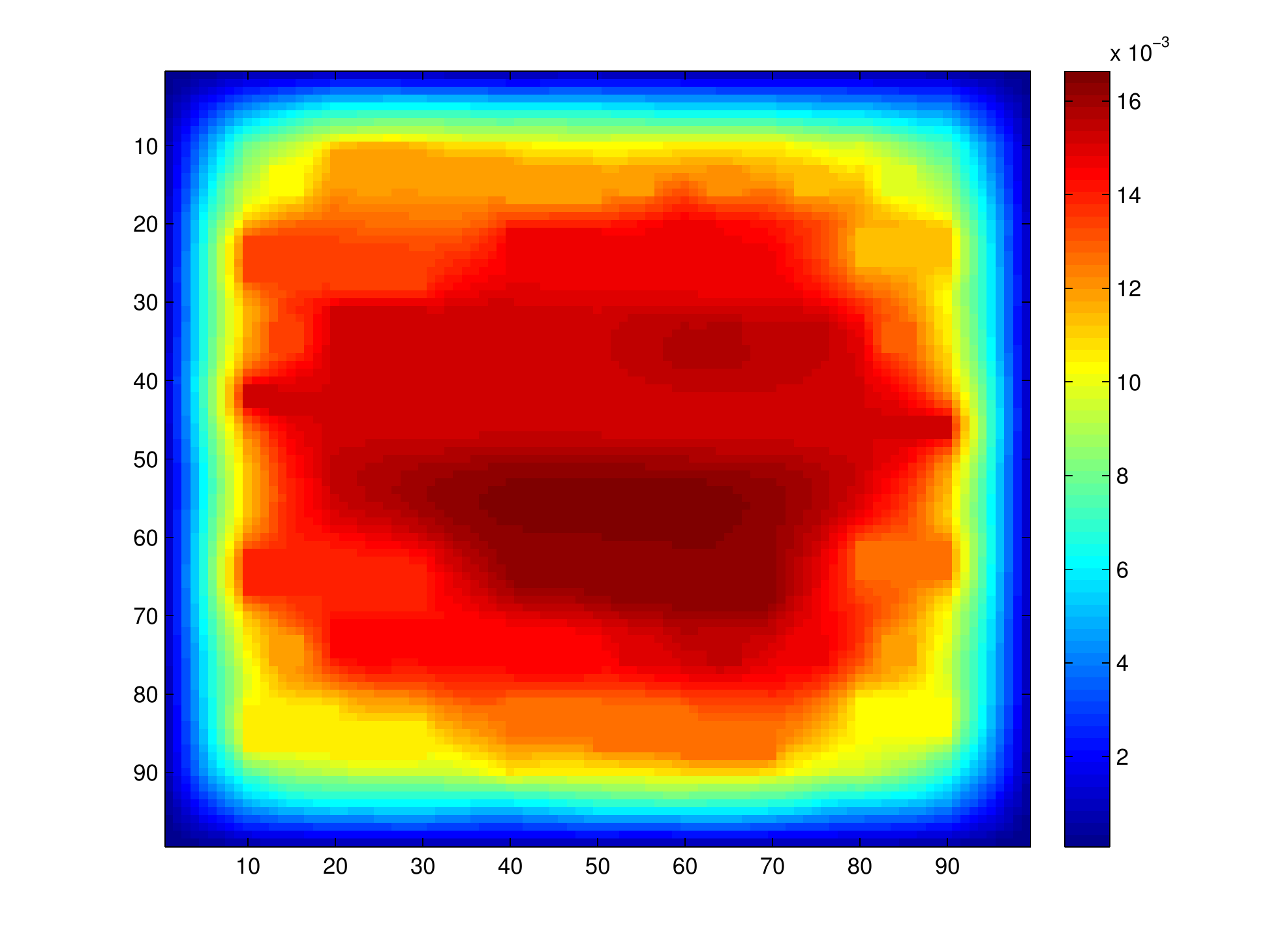}
\includegraphics[width=1.5in, height=1.5in]{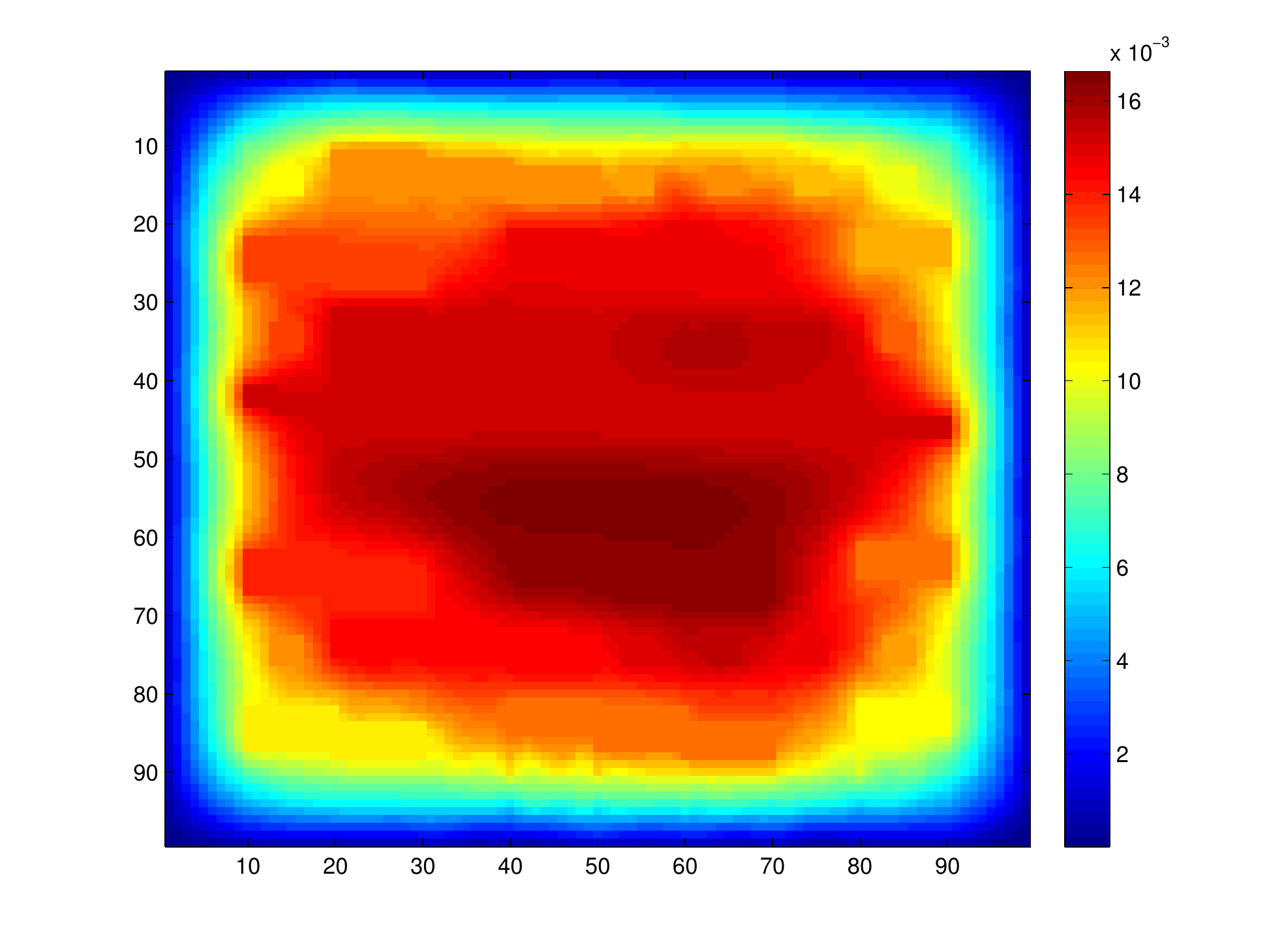}
\caption{Plots of the reference solution (left) and
sample mean of numerical solution at $T = 0.02$: sequential sampling (middle), full sampling (right).}
\label{fig:sol_samp}
\end{figure}

\begin{figure}[ht!]
\centering
\includegraphics[width=0.45\linewidth]{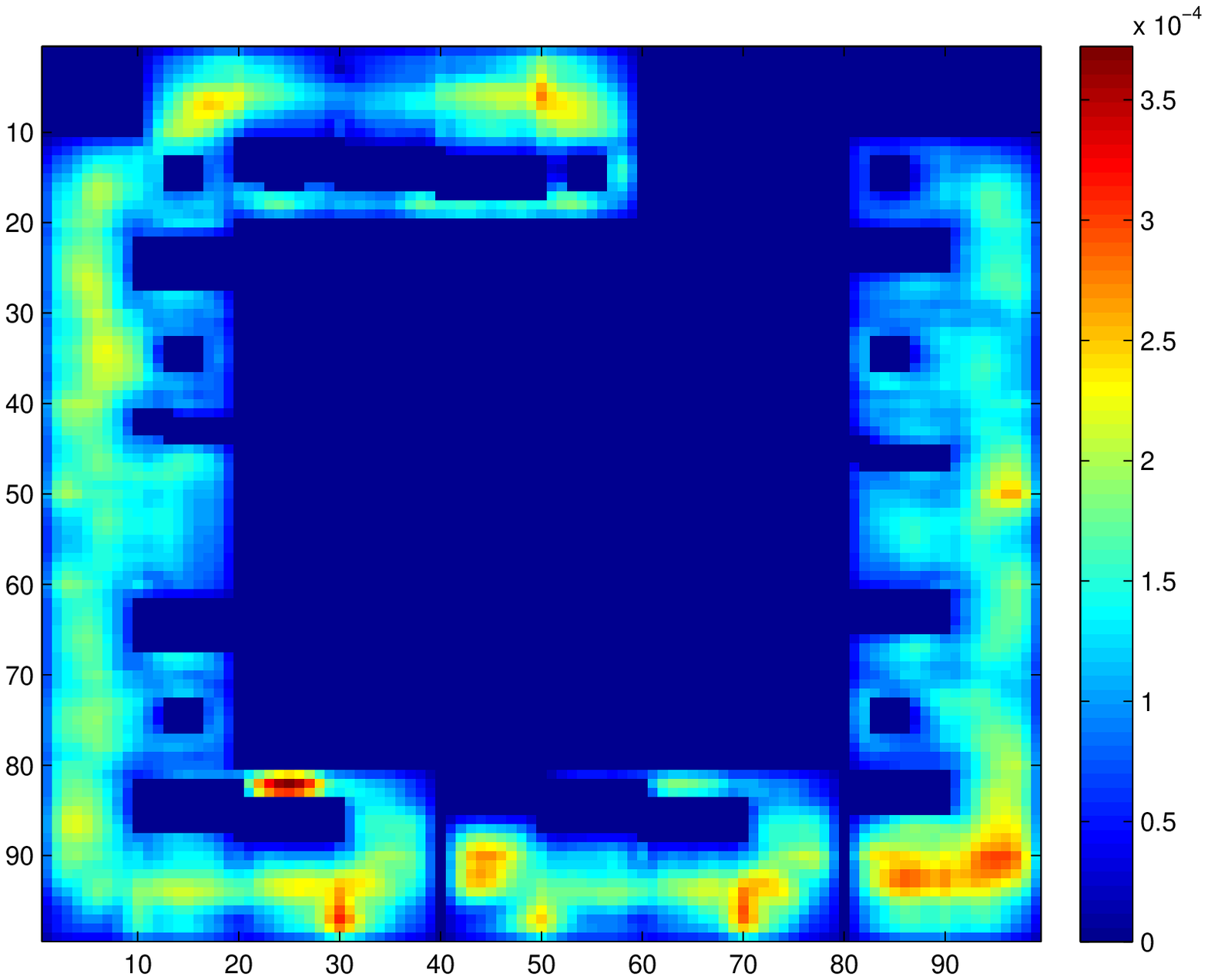}
\includegraphics[width=0.45\linewidth]{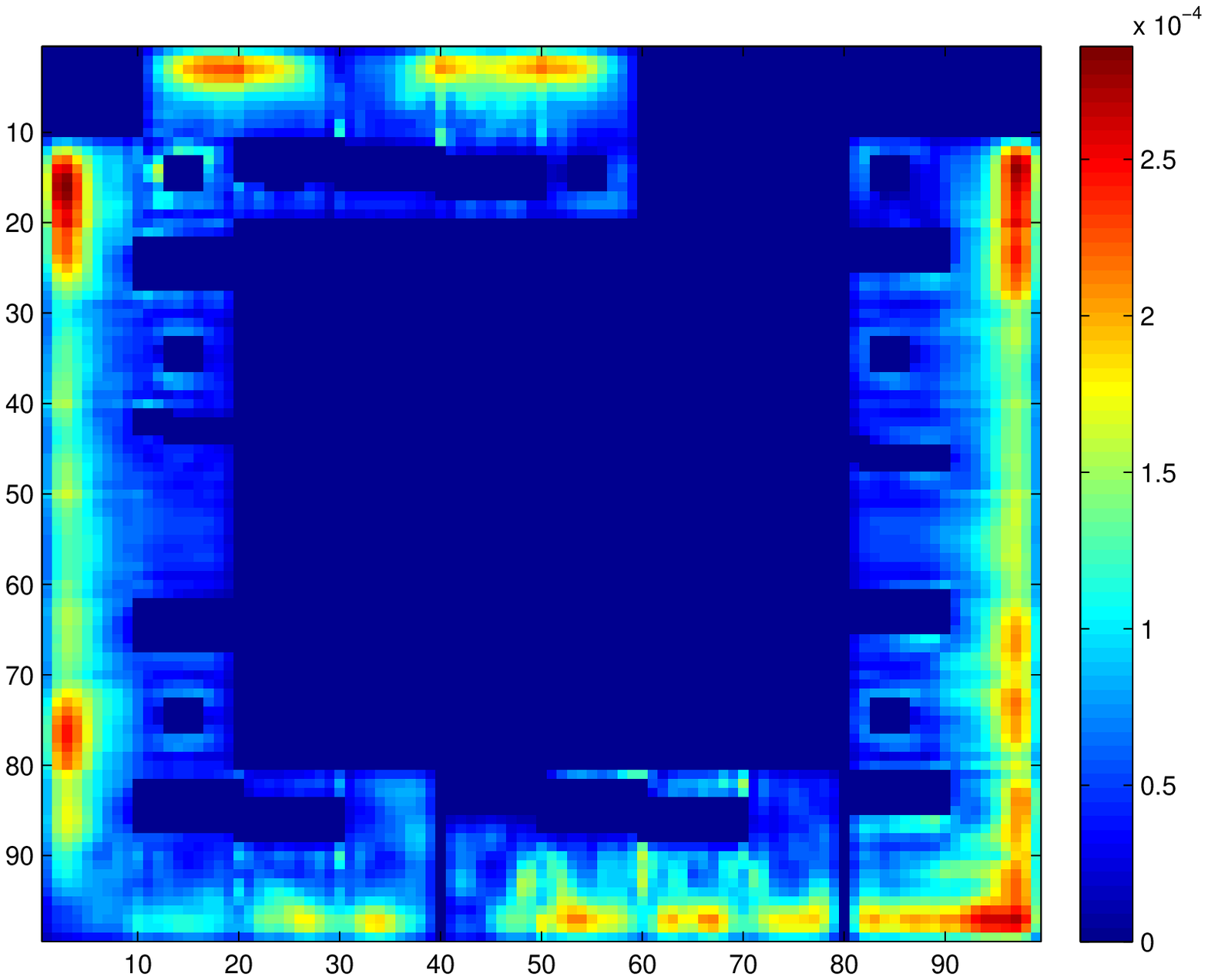}
\caption{Plots of sample standard deviation of numerical solution at $T = 0.02$: sequential sampling (left), full sampling (right).}
\label{fig:sd_samp}
\end{figure}

Figure \ref{fig:sol_samp} shows the reference solution and
the sample mean at $T=0.2$.
The $L^2$ error for the mean at $T = 0.02$ is $0.92\%$
in the full sampling method, lower than $2.24\%$ in the sequential
sampling method.
Figure \ref{fig:sd_samp} shows the pixel-wise standard deviation of
the samples. It can be seen that the deviation is smaller in full sampling.


\begin{figure}[ht!]
\centering
\includegraphics[width=0.45\linewidth]{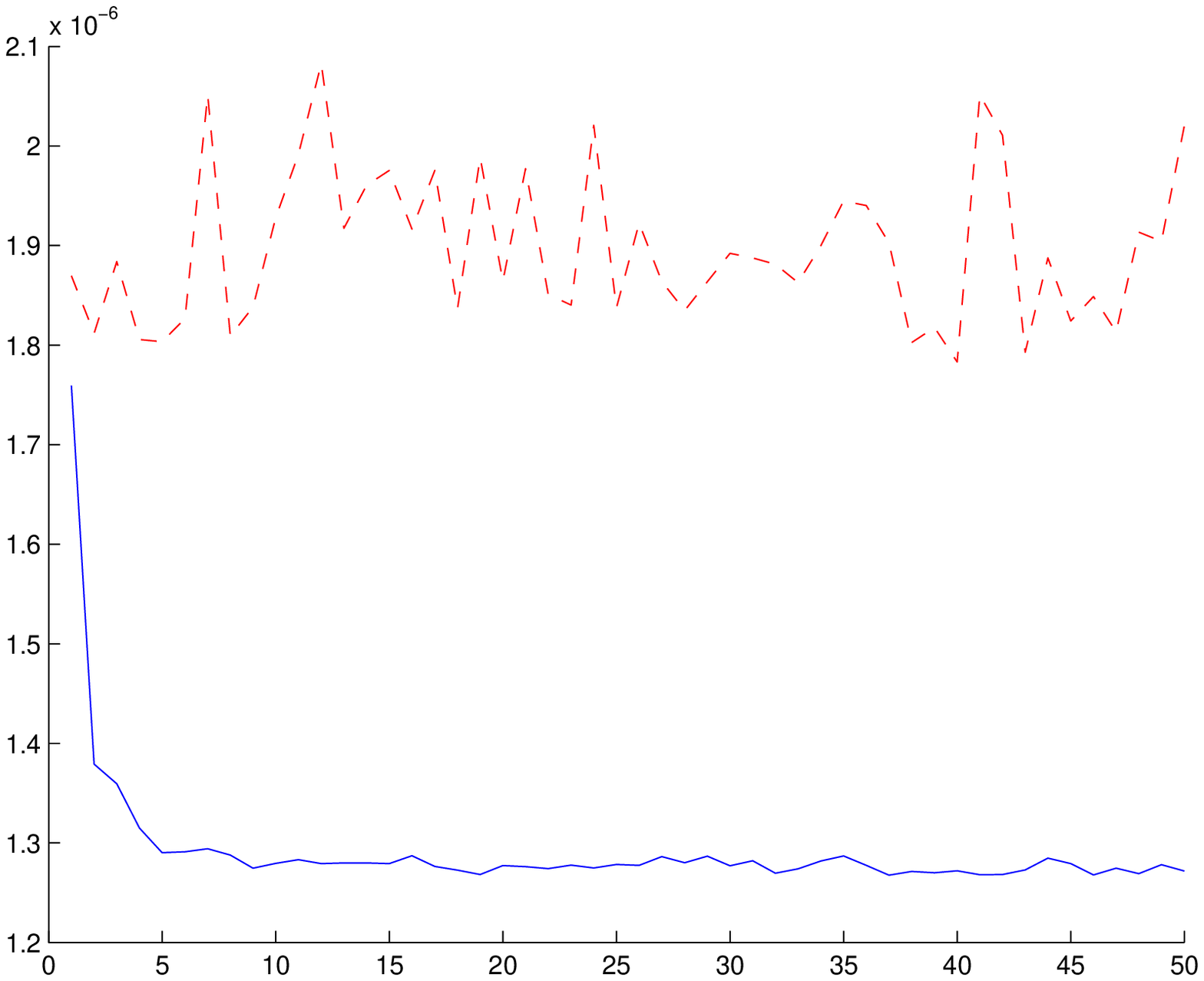}
\includegraphics[width=0.45\linewidth]{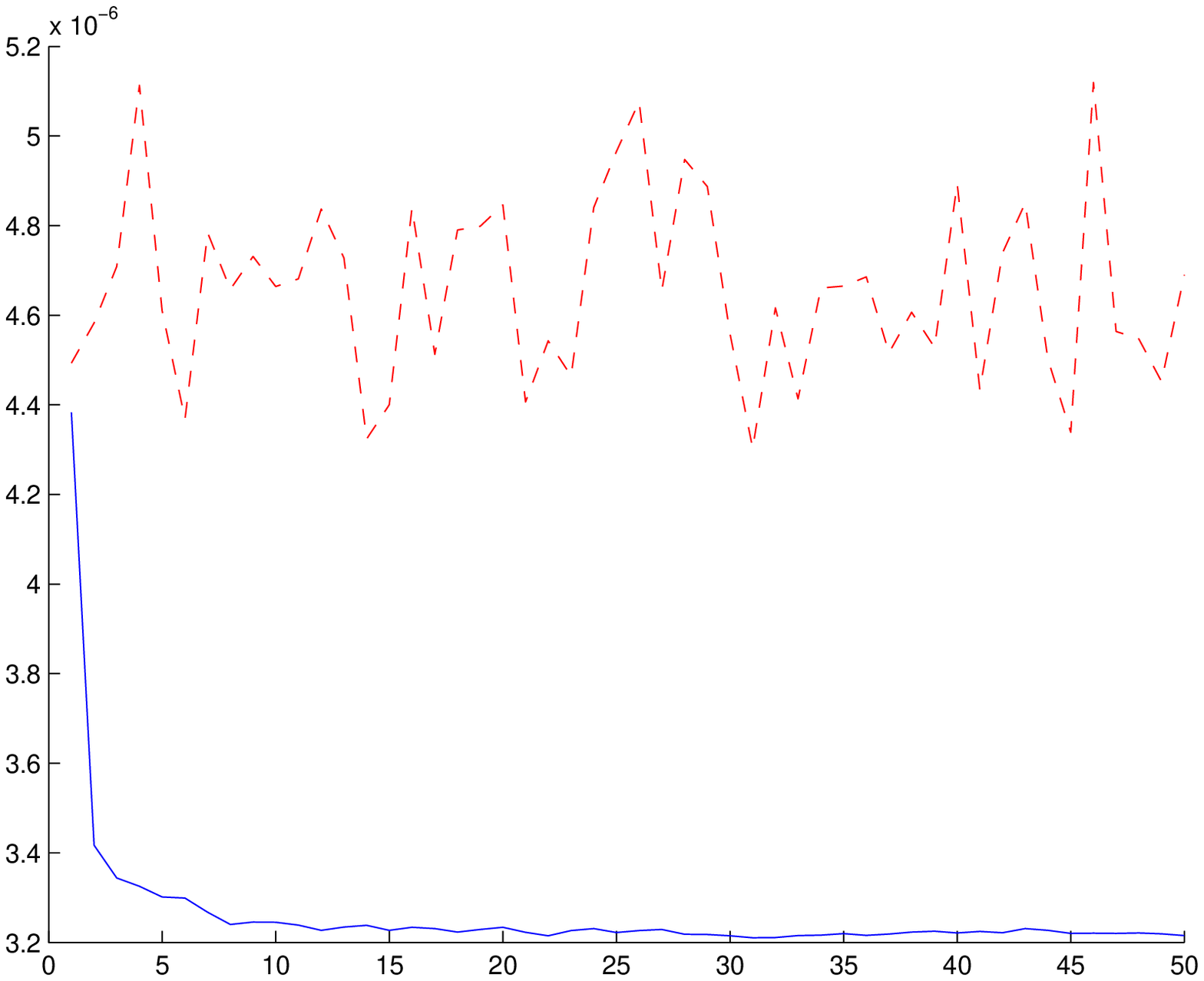}
\caption{Residual vs sample using sequential sampling (red dotted line) and full sampling (blue solid line): at time $T = 0.01$ (left), at time $T = 0.02$ (right).}
\label{fig:res}
\end{figure}

\begin{figure}[ht!]
\centering
\includegraphics[width=0.45\linewidth]{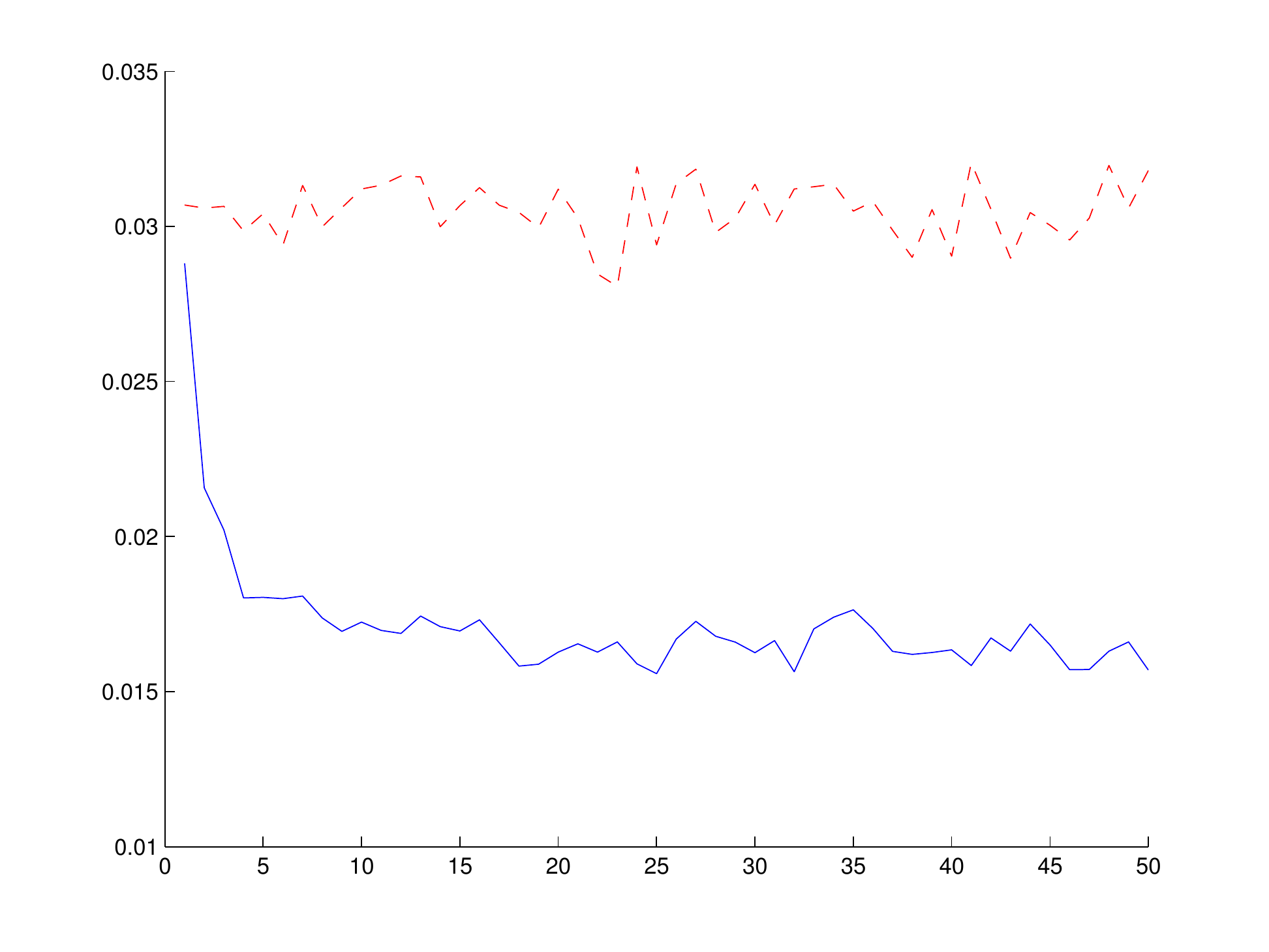}
\includegraphics[width=0.45\linewidth]{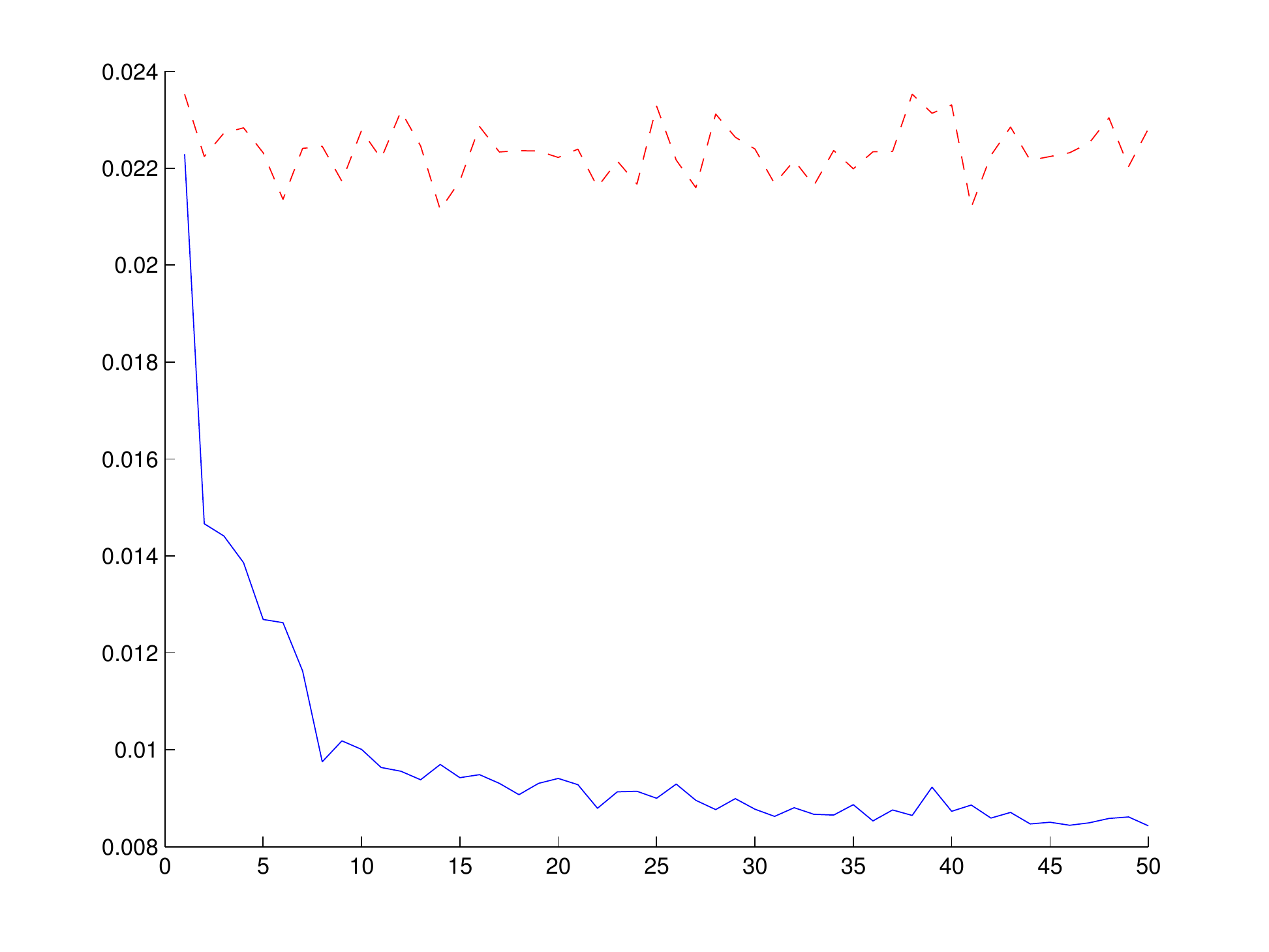}
\caption{$L^2$ error vs sample using sequential sampling (red dotted line) and full sampling (blue solid line): at time $T = 0.01$ (left), at time $T = 0.02$ (right).}
\label{fig:L2}
\end{figure}

\begin{figure}[ht!]
\centering
\includegraphics[width=0.45\linewidth]{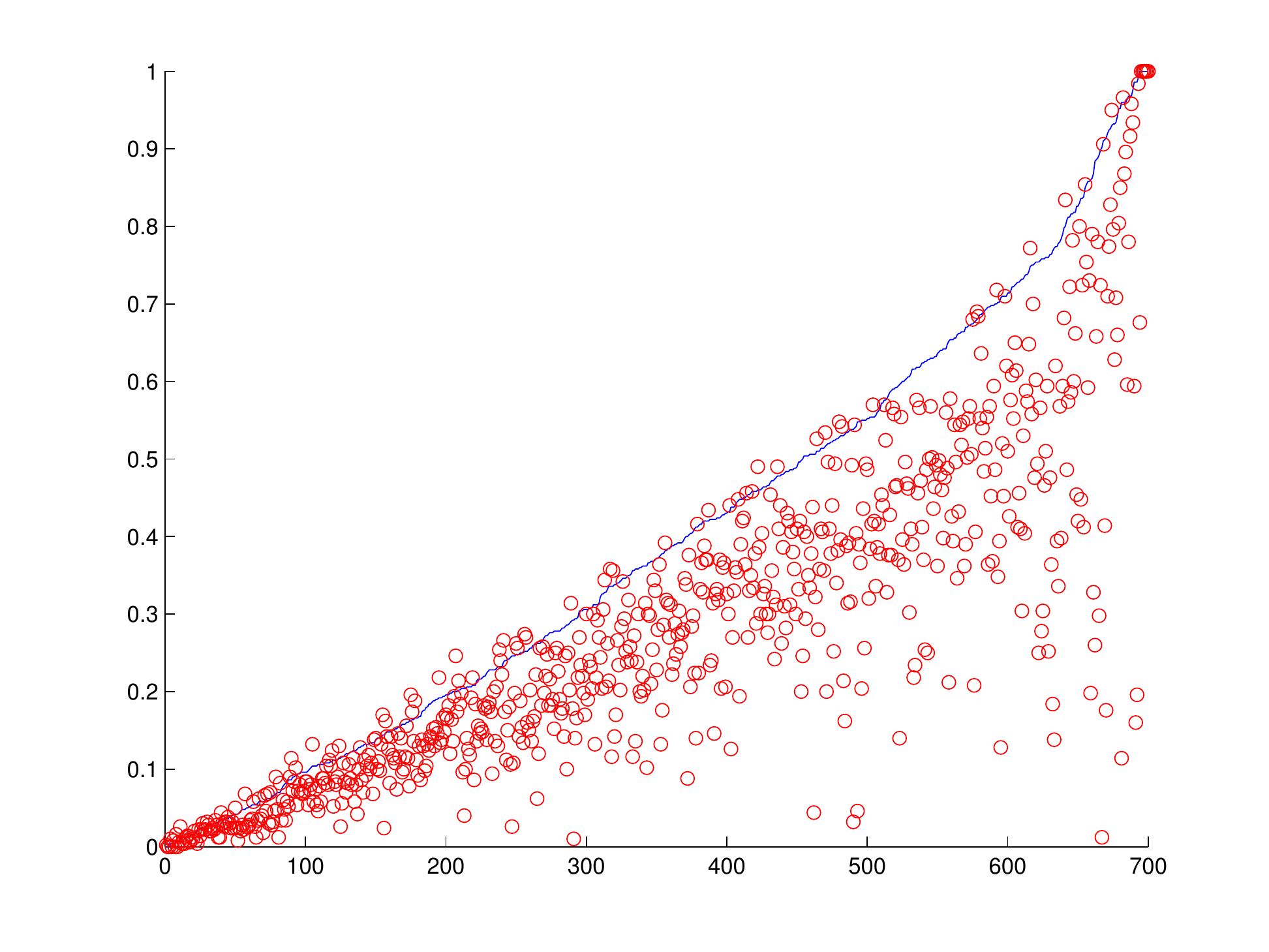}
\includegraphics[width=0.45\linewidth]{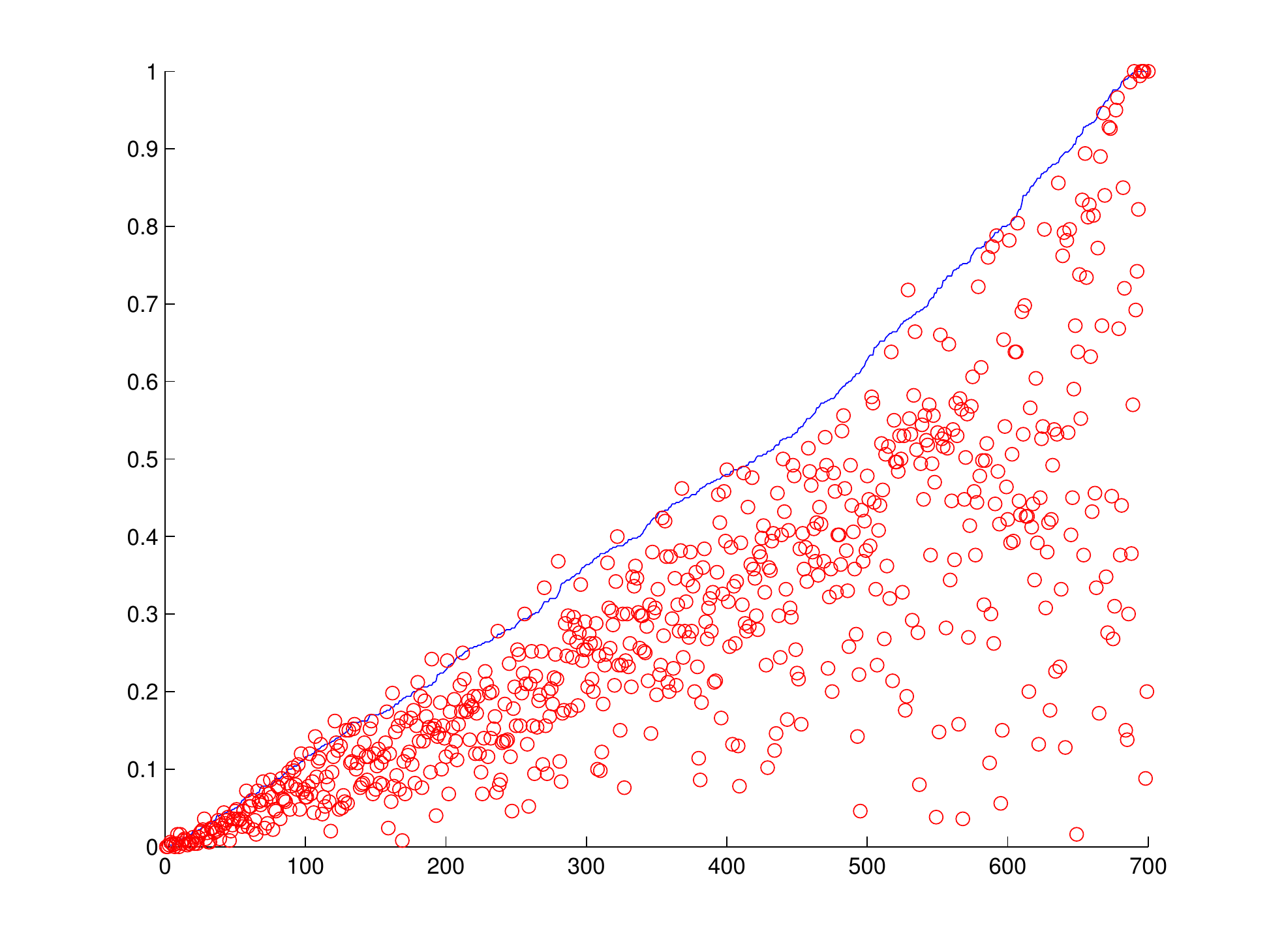}
\caption{History of occurrence probability against basis functions using sequential sampling (red dotted line) and full sampling (blue solid line): at time $T = 0.01$ (left), at time $T = 0.02$ (right).}
\label{fig:dist}
\end{figure}

\begin{figure}[ht!]
\centering
\includegraphics[width=0.45\linewidth]{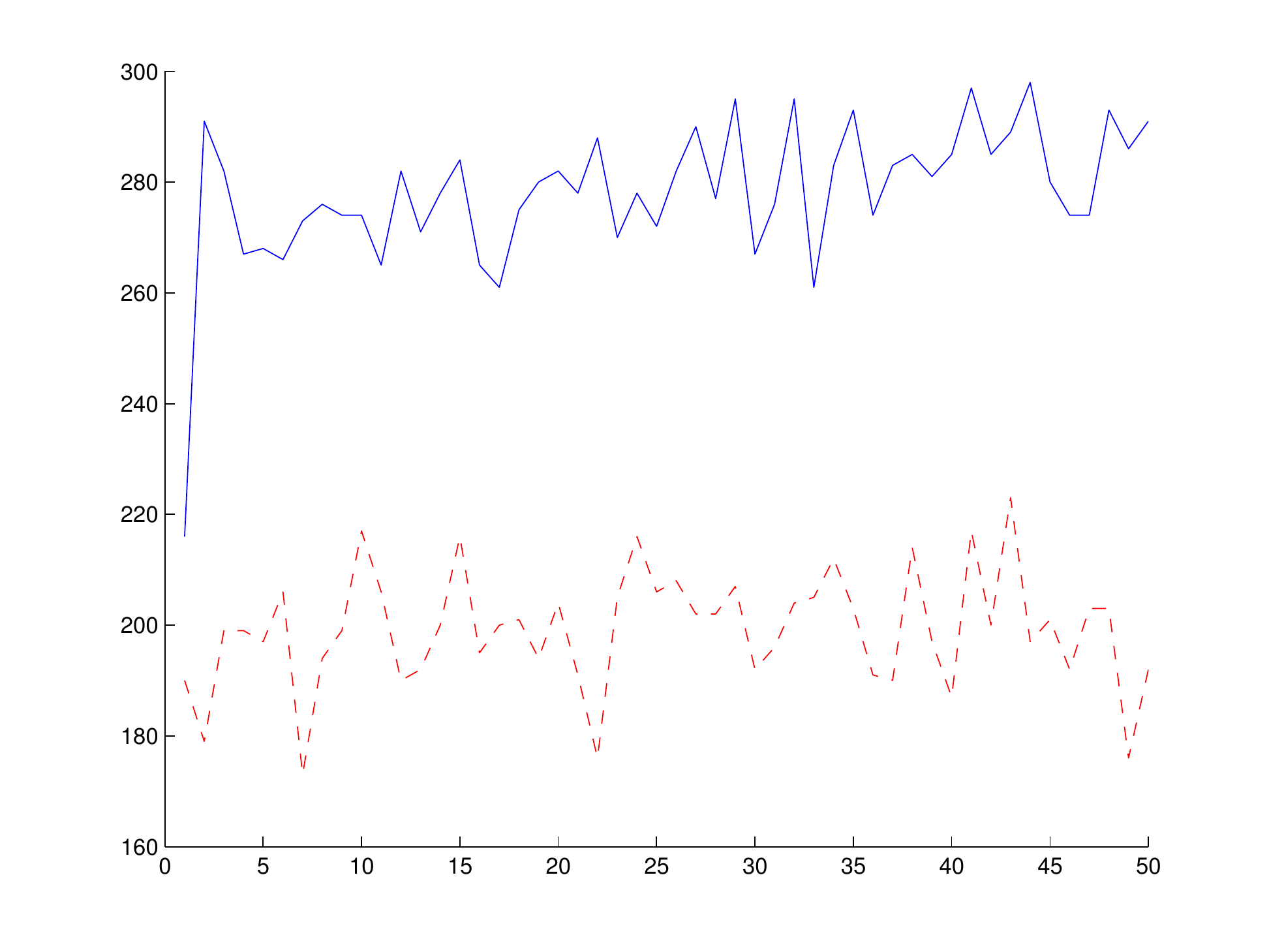}
\includegraphics[width=0.45\linewidth]{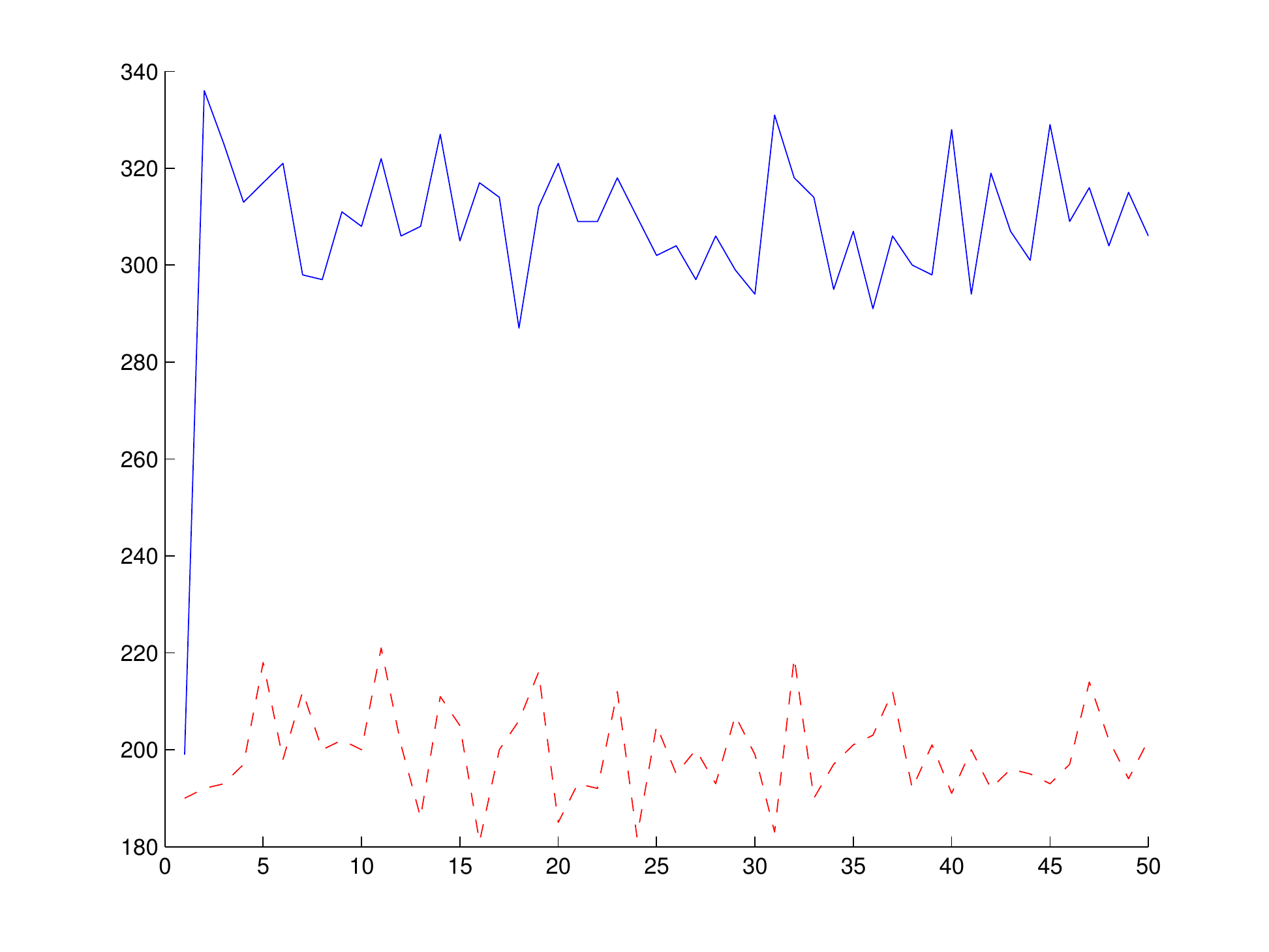}
\caption{History of number of basis functions against sampling process using sequential sampling (red dotted line) and full sampling (blue solid line): at time $T = 0.01$ (left), at time $T = 0.02$ (right).}
\label{fig:Nbasis}
\end{figure}

In Figures \ref{fig:res} and \ref{fig:L2},
the residual and $L^2$ errors are plotted for
both sequential and full sampling. We observe that
the errors and the residual in full sampling decrease
and stabilize in a few iterations.
Moreover, the full sampling gives more accurate solutions
associated with our error threshold in the residual.
In Figure \ref{fig:dist}, we plot the frequency of the basis functions
vs. sample
that appear in full and sequential sampling. The correlation
of the frequencies is $0.87$, which indicates that we have a
good prior model based on the residual.
Finally, in Figure \ref{fig:Nbasis}, we show the full number of basis
functions. As we observe that the number of basis functions in full
sampling approaches to a steady state, which depends on the error
threshold $\sigma_L$ (see (\ref{eq:post_sigma})).


\subsection{Wave equation}
\label{sec:numerics_wave}

In this section, we present an application of
our approach to wave equations.
Multiscale basis function construction follows to a similar procedure
described in Appendix \ref{sec:app1} with some slight modifications,
see \cite{chung2014generalized} for details. Here,
we describe the residual that is used
in our Bayesian framework.

We will use IPDG method with multiscale basis functions to solve the following wave equation,
\begin{equation}
\frac{\partial^{2}u}{\partial t^{2}}=\nabla\cdot(a\nabla u)+f\quad\quad\mbox{ in }\quad \Omega \times [0,T].
\label{eq:waveeqn}\end{equation}
The IPDG method is to find $u_{H}$ such that
 \begin{equation}
\int_\Omega \frac{\partial^2 u_{H}}{\partial t^{2}} \, v+a_{DG}(u_{H},v)=\int_\Omega f\, v ,\quad\quad\forall\, v\in V_{H},\label{eq:ipdg}
\end{equation}
where the bilinear form $a_{DG}(u,v)$ is defined by
\begin{eqnarray*}
a_{DG}(u,v) & =& \sum_{K\in\mathcal{T}^{H}}\int_{K}a\nabla u\cdot\nabla v+\sum_{e\in\mathcal{E}^{H}}(-\int_{e}\{a\nabla u\cdot n\}_{e}\,[v]_{e} \\
 & & -\int_{e}\{a\nabla v\cdot n\}_{e}\,[u]_{e}+\cfrac{\gamma}{h}\int_{e}a[u]_{e}\,[v]_{e}),
\end{eqnarray*}
 with $\gamma>0$ is a penalty parameter
 and $n$ denotes the unit
normal vector on $e$.
Let $\Delta t > 0$ be the time step size.
We will consider  the classical second order central finite difference method for time discretization, i.e., we find $u^{n+1}_H$
such that
\begin{equation*}
\int_\Omega\frac{ u^{n+1}_{H} - 2u^n_{H} + u^{n-1}_{H} }{\Delta t^2}v + a_{DG}(u^n_{H},v) = \int_\Omega f^n\, v \quad\quad\forall\, v\in V_{H}.
\label{eq:fulldiscrete}
\end{equation*}

We compute the solution $u_{H}^{n+1}$ at $n+1$ th time level, starting with the offline space.
In each coarse element $K\in\mathcal{T}^{H}$, we define the local
residual operator, $R_{v}^{n}$, as
\[
R_{v}^{n}(u^{n+1}_H,u^n_H,u^{n-1}_H)=(\cfrac{u_{H}^{n+1}-2u_{H}^{n}+u_{H}^{n-1}}{\Delta t^{2}},v)+a_{DG}(u_{H}^{n},v)-(f^n,v),
\]
and the operator norm of local residual functional is defined by
\[
\|\mathcal{R}_{K}^{n}\|=\sup_{v\in V_{h}(K)}\cfrac{|R_{v}^{n}|}{\|v\|_{L^{2}(K)}}.
\]

\begin{figure}[H]
\begin{centering}
\includegraphics[scale=0.5]{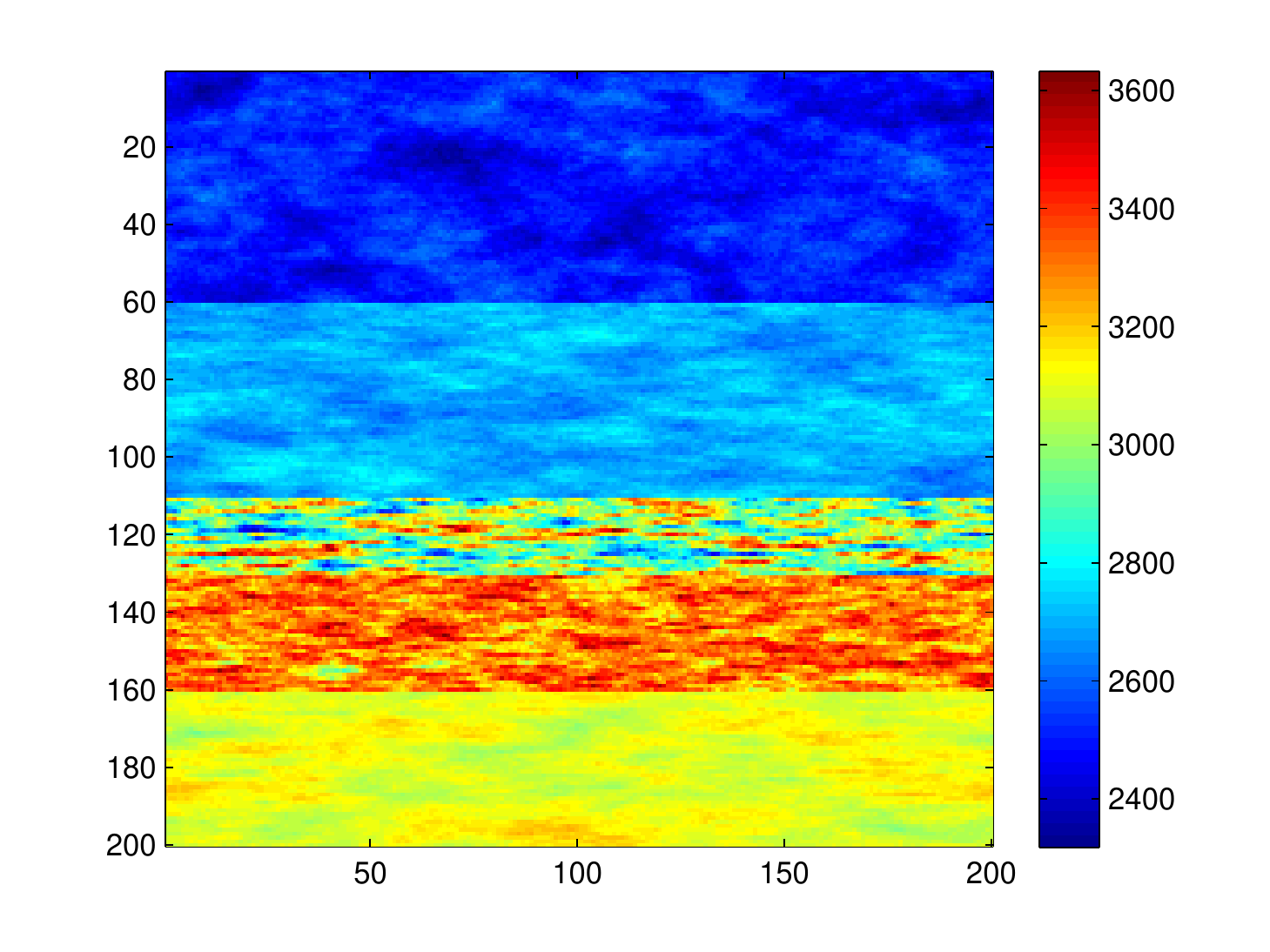}
\par\end{centering}

\protect\caption{The medium $a$.}
\label{fig:medium_wave}
\end{figure}

We use the medium
$a$ shown in Figure \ref{fig:medium_wave}.
In this example, we find the distribution of the solution
at the time instants $T=0.4$.
The fine grid is $200\times 200$ and the coarse grid is $20\times 20$.
We use one permanent basis function per edge to
compute ``fixed'' solution and use our Bayesian framework
to seek additional basis functions by solving small global
problems.
In this approach, we use previous time solution in the posterior.
In our approach, using the global residual, some local regions
($17$ \% of total coarse regions)
are defined,
where multiscale basis functions are added.
 In these coarse blocks,
we apply both sequential sampling and full sampling algorithms.

\begin{figure}[H]
\begin{centering}
\includegraphics[scale=0.25]{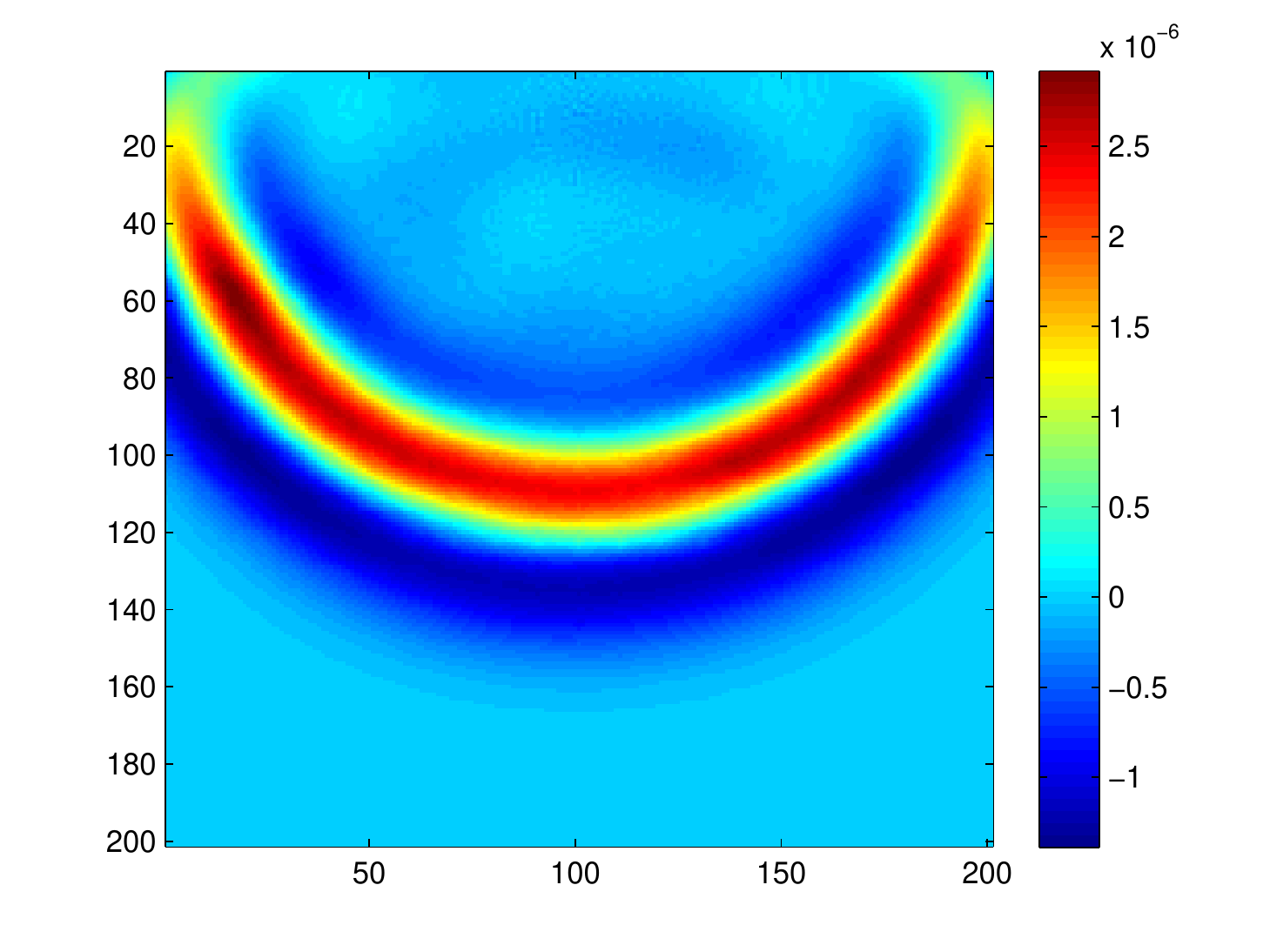} \includegraphics[scale=0.25]{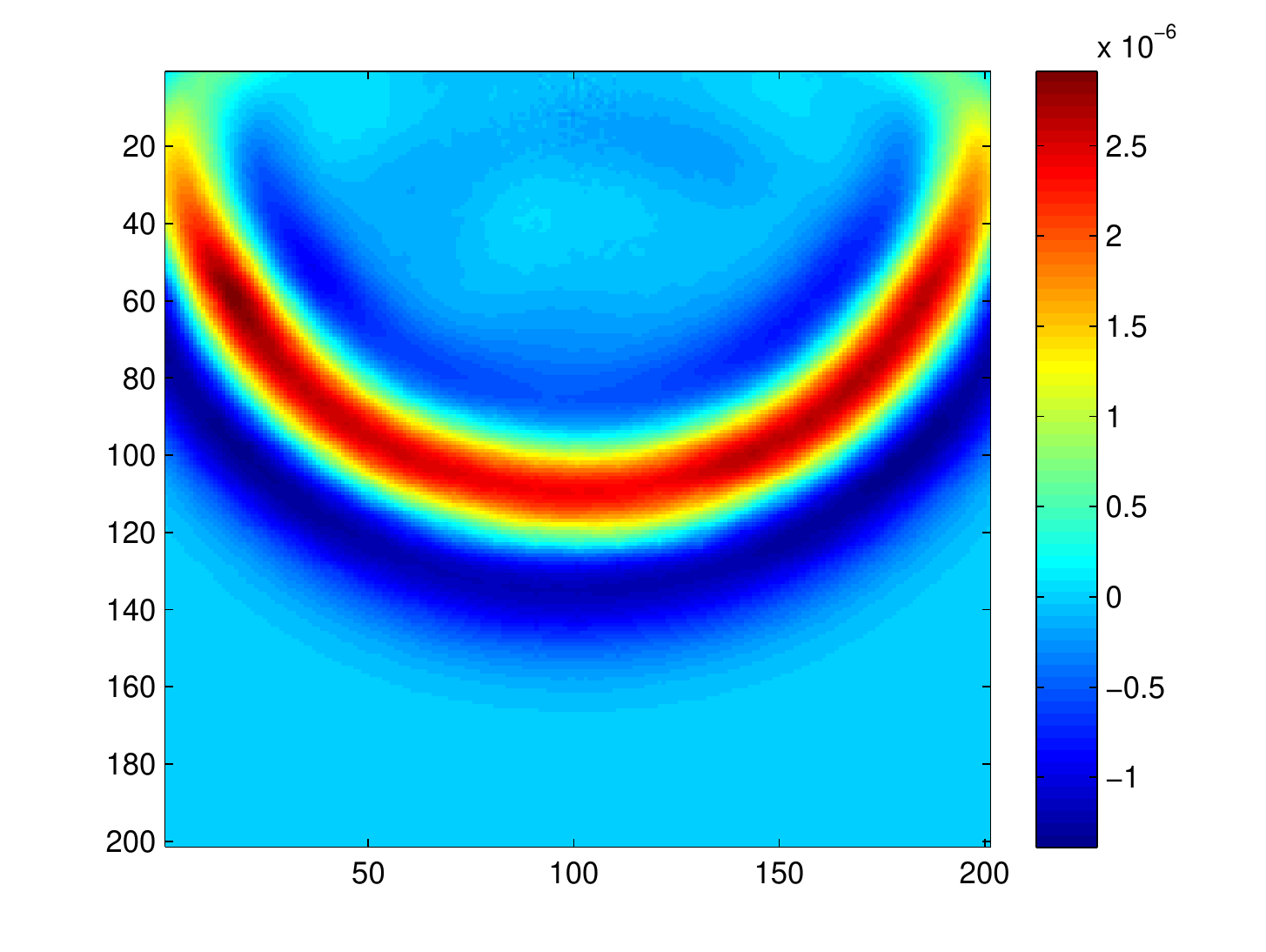} \includegraphics[scale=0.25]{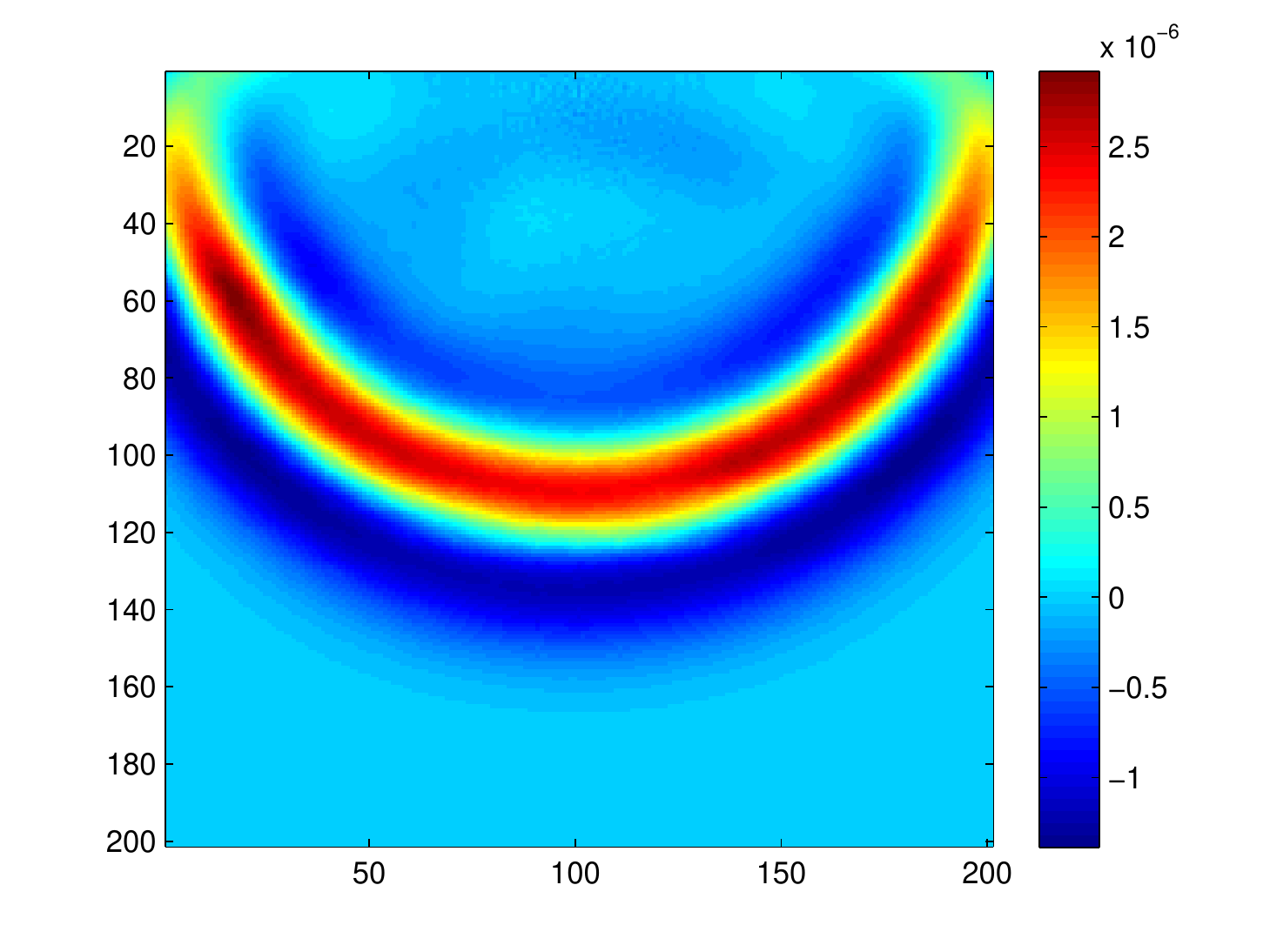}
\par\end{centering}
\protect\caption{Numerical solution at $T=0.4$. Left: reference solution, Middle: mean
solution for sequential sampling, Right: mean solution for full sampling.}
\label{fig:mean_alg1_Ex2}
\end{figure}

\begin{figure}[H]
\begin{centering}
\includegraphics[scale=0.4]{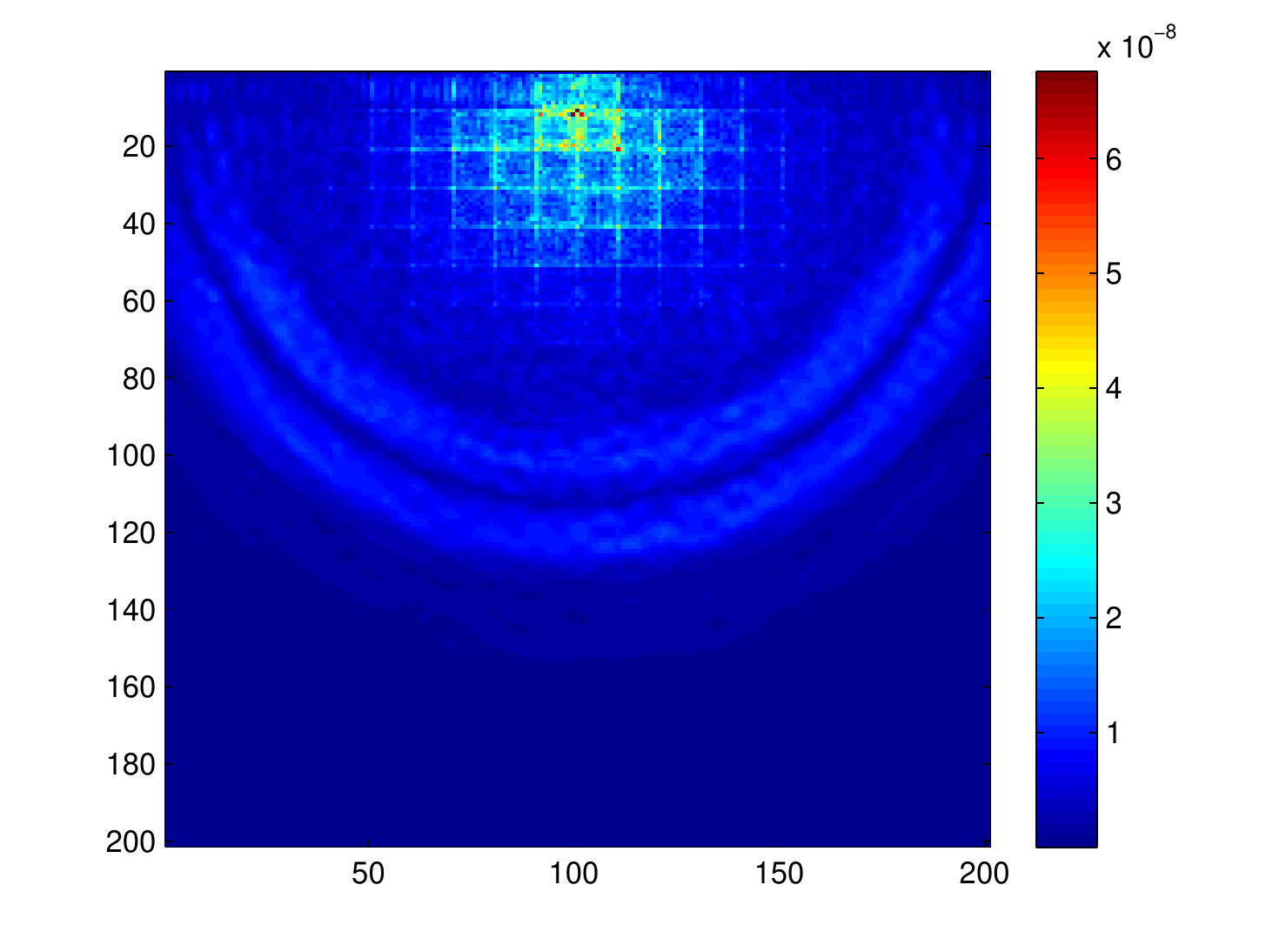}\includegraphics[scale=0.4]{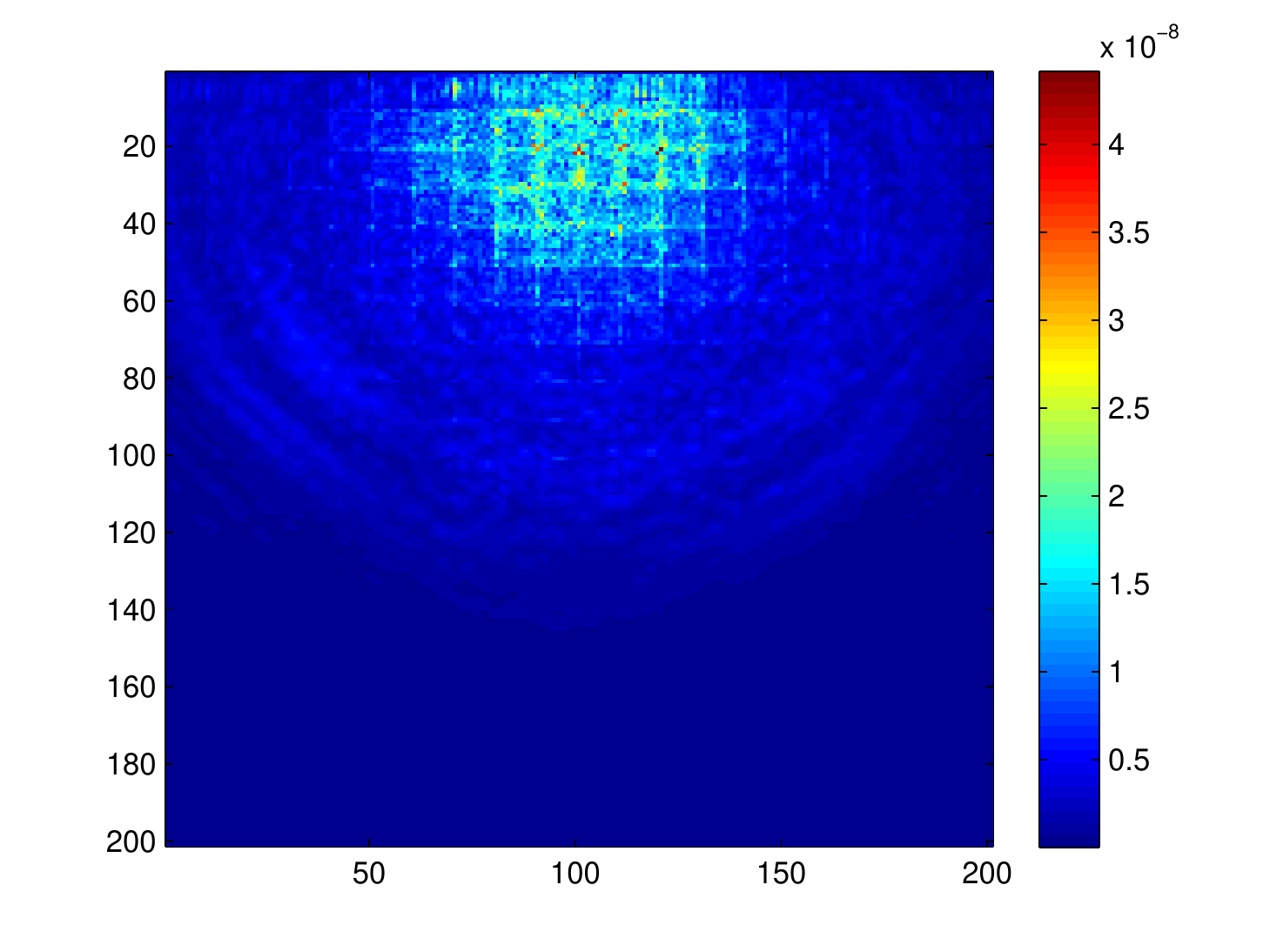}
\par\end{centering}

\protect\caption{Standard deviation of the solution. Left: sequential sampling, Right:
full sampling}
\label{fig:std_alg1_Ex2}
\end{figure}

\begin{figure}[H]
\begin{centering}
\includegraphics[scale=0.4]{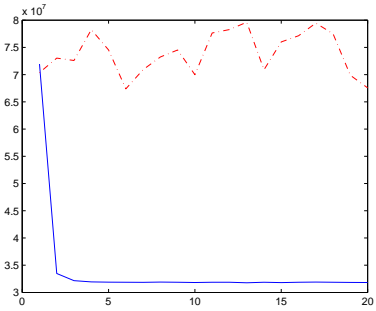}
\par\end{centering}

\protect\caption{Residual vs sample using sequential sampling (red dotted line) and full sampling (blue solid line) at the final time.}
\label{fig:residual_Ex2}
\end{figure}

\begin{figure}[H]
\begin{centering}
\includegraphics[scale=0.4]{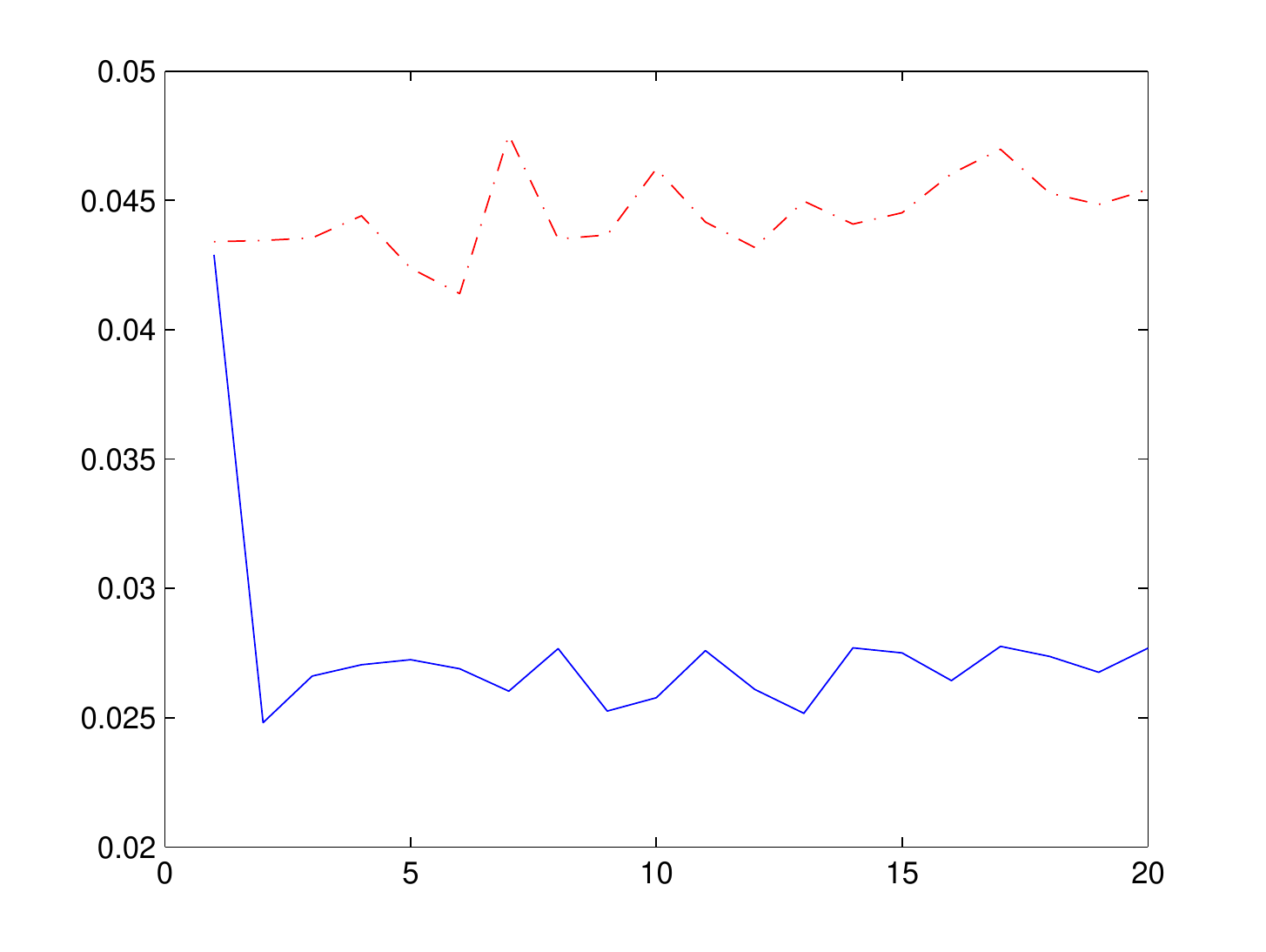}
\par\end{centering}

\protect\caption{$L^2$ error vs sample using sequential sampling (red dotted line) and full sampling (blue solid line) at the final time.}
\label{fig:error_Ex2}
\end{figure}

\begin{figure}[H]
\begin{centering}
\includegraphics[scale=0.4]{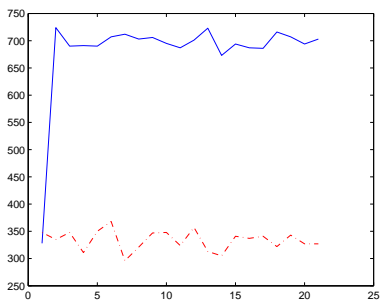}
\par\end{centering}
\protect\caption{History of number of basis functions against sampling process using sequential sampling (red dotted line) and full sampling (blue solid line): at time T = 0.4.}
\label{fig:num_basis_Ex2}
\end{figure}

In Figure \ref{fig:mean_alg1_Ex2}, we depict the mean
solution using the sequential sampling algorithm and full sampling
algorithms. The errors for the mean are $4.47\%$ for the sequential sampling
and $2.67\%$ for full sampling.
As before, the full sampling algorithm provides a more accurate result.
We depict the standard deviation of the solution at each pixel
in Figure \ref{fig:std_alg1_Ex2}. We observe that the
true solution falls within the limits of the
 mean and the standard
deviations. Next, we show the results across several samples.
For the sequential sampling, we use $20$ realizations and show both
residuals and the errors in Figure \ref{fig:residual_Ex2} and
\ref{fig:error_Ex2}. The error is computed
as a difference between the solution
and the snapshot solution using the snapshot vectors in the elements,
where the basis functions are updated.
From these figures, we observe that the residuals and errors are
smaller for full sampling compared to sequential sampling. Moreover,
Gibbs sampling stabilizes in a few iterations.
In Figure \ref{fig:num_basis_Ex2}, we plot the total number of
basis functions for samples. As we observe that the full sampling
results to about $700$ additional basis functions (on average)
for the entire domain. This number is affected by the error threshold
$\sigma_L$ (see (\ref{eq:post_sigma}))
that we impose in the posterior distribution.

\section{Conclusions}

In this paper, we propose a novel Bayesian Multiscale Finite Element
approach. The main idea of the approach is to solve time-dependent
problems, e.g., parabolic equations or
wave equations, using multiscale basis functions in each coarse grid.
It is known that \cite{chung2016adaptive}
 first few basis functions reduce the error
substantially,
while the ``rest'' of basis functions may reduce error significantly
less. This is due to global effects. Various online approaches
(e.g., \cite{chung2016adaptive} and the references therein)
are proposed to remedy this situation. The online approaches
are expensive as we need to compute new basis functions.
In this paper, we propose a rigorous sampling for ``rest'' of
basis functions, which allow computing multiple realizations of
the solution and thus provide a probabilistic description for
un-resolved scales. Our approaches are motivated by recent
work \cite{mallick2016}, where the authors propose a probabilistic
numerical method.

In the paper, we present a description of the method
and sampling algorithms for un-resolved scales.
Several posterior distributions are considered and two
sampling mechanisms are studied. In both sampling methods,
we use local residuals as a prior distributions for
selecting multiscale basis functions. In the first approach,
multiscale basis functions are selected from the prior distribution,
while in the second approach, we present a full sampling using
the Gibbs sampling. We show that the Gibbs sampling rapidly stabilizes
at the steady state and provides a more accurate approximation
at an additional cost.
We consider several discretizations and applications to wave and parabolic
equations. In general, the proposed method can be used to
condition the solutions at measurement locations or
computing solutions in stochastic environments.
In this case, one can combine the precisions of the forward and inverse
simulations in a unified way.
This will be studied
in our future work.

\appendix\section{Some Notations}

\begin{itemize}

\item $\omega_i$ is a subdomain with a common vertex at $x_i$;
$\omega_E$ is a subdomain with a common edge $E$;
$K$ is a coarse-grid block

\item $u_H$ is the coarse-grid multiscale solution;
$u_H^{n}$ is the coarse-grid space-time multiscale solution defined
on $(T_{n-1},T_n)$

\item $\psi_i^{\omega_j}$  - local snapshot solutions in
the oversampled region $\omega_j$; $\phi_i^{\omega_j}$  - local offline basis in $\omega_j$

\item $\phi_i^{n,\omega_j}$  - permanent (first few)
local offline basis in $\omega_j$;
$\phi_{i,+}^{n,\omega_j}$  - the rest of basis functions in $\omega_j$

\item $u_H^{n,\text{fix}}$ is the solution in
$(T_{n-1},T_n)$ computed using permanent basis functions

\item
$R^{n,\text{fix}}_v(u_{H}^{n,\text{fix}}(x,t),u_{H}^{n,\text{fix}}(x,T_{n-1}))$ -
the residual corresponding to fixed
solution, $v$ is the test function

\item $\mathcal{R}^n$ - the global residual vector;
$\mathcal{R}^n_{\omega_j}$ - the local residual vector  in $\omega_j$

\item $V_{H,\text{off}}$ is the offline space in $\Omega$;
$V_{H,\text{off}}^{\omega_i}$ is the offline space in $\omega_i$;
$V_{H,\text{snap}}$ is the snapshot space in $\Omega$;
$V_{H,\text{snap}}^{\omega_i}$ is the snapshot space in $\omega_i$

\item $V_{H,\text{off}}^{(T_1,T_2)}$ is the offline space in $\Omega \times (T_1,T_2)$;
$V_{H,\text{off}}^{\omega_i,(T_1,T_2)}$ is the offline space in $\omega_i \times (T_1,T_2)$; $V_{H,\text{snap}}^{(T_1,T_2)}$ is the snapshot space in $\Omega \times (T_1,T_2)$;
$V_{H,\text{snap}}^{\omega_i,(T_1,T_2)}$ is the snapshot space in $\omega_i \times (T_1,T_2)$

\end{itemize}

\appendix\section{Space-time basis construction details}
\label{sec:app1}

\subsection{Snapshot space}

Let $\omega$ be a coarse neighborhood in space.
Coarse node index is omitted to simplify the notations.
The construction of the
offline basis functions in $(T_{n-1},T_n)$ uses
a snapshot space $V_{H,\text{snap}}^{\omega}$
(or $V_{H,\text{snap}}^{\omega,(T_{n-1},T_n)}$).
For simplicity,
the coarse time index $n$ is omitted.
The snapshot space
$V_{H,\text{snap}}^{\omega}$ consists of functions in $\omega$
and contains all or most necessary components of the fine-scale
solution restricted to $\omega$. Further, a spectral problem is
solved
in the snapshot space to compute
multiscale basis functions (offline space).

We consider a snapshot space that consists
of solving local problems for all possible boundary conditions.
We denote by $\omega^{+}$ the oversampled space region of
$\omega \subset\omega^{+}$, defined by adding several fine- or coarse-grid
layers around $\omega$ and define $(T_{n-1}^{*}, T_{n})$ as
the left-side oversampled time region for $(T_{n-1},T_{n})$.
We can compute inexpensive snapshots using random boundary conditions on
the oversampled space-time region $\omega^{+}\times(T_{n-1}^{*},T_{n})$.
by solving
a small number of local problems imposed with random boundary conditions
\begin{equation*}
\begin{split}
-\text{div} (\kappa(x,t_*) \nabla \psi_{j}^{+,\omega})=0\ \ \text{in}\ \omega^{+}\times (T_{n-1}^{*},T_{n}), \\
\psi_{j}^{+,\omega}(x,t)= r_l\  \ \text{on} \ \ \partial \left( \omega^{+}\times (T_{n-1}^{*},T_{n}) \right),
\end{split}
\end{equation*}
where $t_*$ is a time instant,
$r_l$ are independent identically distributed (i.i.d.) standard Gaussian random vectors on the fine-grid nodes of the boundaries on $\partial \omega^{+}\times (T_{n-1}^{*},T_{n})$.
Then the local snapshot space on $\omega^{+}\times (T_{n-1}^{*},T_{n})$ is
\[
V_{H,\text{snap}}^{+,\omega} = \text{span}\{\psi_{j}^{+,\omega}(x,t) | j=1,\cdot\cdot\cdot, l^{\omega}+p_{\text{bf}}^{\omega}\},
\]
where $l^{\omega}$ is the number of local offline basis we want to construct in $\omega$ and $p_{\text{bf}}^{\omega}$ is the buffer number.

\subsection{Offline space}

We perform a space reduction by appropriate spectral problems to compute
the offline space.
We solve $(\phi,\lambda)\in V_{H,\text{snap}}^{+,\omega}\times\mathbb{R}$ such that
\begin{equation}\label{eq:eig-problem}
A_n(\phi,v) = \lambda S_n(\phi,v), \quad \forall v \in V_{\text{snap}}^{\omega^{+}},
\end{equation}
where the bilinear operators $A_n(\phi,v)$ and $S_n(\phi,v)$ are defined by
\begin{equation}
\begin{split}
A_n(\phi,v) =
 \int_{T_{n-1}}^{T_{n}}\int_{\omega^{+}}\kappa(x,t_*)\nabla\phi \cdot \nabla v, \\
S_n(\phi,v) =\int_{T_{n-1}}^{T_{n}}\int_{\omega^{+}}\widetilde{\kappa}^{+}(x,t_*)\phi v,
\end{split}
\end{equation}
where the weight function $\widetilde{\kappa}^{+}(x,t)$ is defined by
$\widetilde{\kappa}^{+}(x,t) = \kappa(x,t)\sum_{i=1}^{N_c}|\nabla\chi_i^{+}|^2$,
$\{\chi_i^{+}\}_{i=1}^{N_c}$ is a partition of unity associated with the oversampled coarse neighborhoods $\{\omega_i^{+}\}_{i=1}^{N_c}$ and satisfies $|\nabla\chi_i^{+}|\geq|\nabla\chi_i|$ on $\omega_i$,
where $\chi_i$ is the standard multiscale basis function for the coarse node $x_i$,
$-\text{div}(\kappa(x,T_{n-1})\nabla\chi_i) = 0$, in $K$,
 $\chi_i=g_i$ on $\partial K$,
for all $K\in\omega_i$, where $g_i$
is linear on each edge of $\partial K$.

We arrange the eigenvalues $\{\lambda_j^{\omega}|j=1,2,\cdot\cdot\cdot\,L_{\omega}+p_{\text{bf}}^{\omega}\}$ from (\ref{eq:eig-problem}) in the ascending order, and select the first $L_{\omega}$ eigenfunctions, which are corresponding to the first $L_{\omega}$ ordered eigenvalues, and
then we can obtain the dominant modes $\psi_{j}^{\omega}(x,t)$ on the target
region $\omega\times (T_{n-1},T_{n})$ by restricting
$\psi_j^{+,\omega}(x,t)$ onto $\omega\times (T_{n-1},T_{n})$. Finally, the offline basis functions on $\omega\times (T_{n-1},T_{n})$ are defined by $\phi_j^{\omega}(x,t) = \chi^\omega\psi_j^{\omega}(x,t)$, where $\chi^\omega$ is the standard multiscale basis function  for a generic coarse neighborhood $\omega$. This product gives conforming basis functions (Discontinuous Galerkin discretizations can also be used). We also define the local offline space on $\omega\times (T_{n-1},T_{n})$ as
\[
V_{H,\text{off}}^{\omega} = \text{span}\{\phi_{j}^{\omega}(x,t) | j=1,\cdot\cdot\cdot, l^{\omega} \}.
\]

Note that one can take $V_{H,\text{off}}^{(T_{n-1},T_n)}$ in the coarse-grid
equation
as $V_{H,\text{off}}^{(T_{n-1},T_n)} =  \text{span}\{\phi_{j}^{\omega_i}(x,t) | 1\leq i\leq N_c, 1\leq j\leq l_{i} \}$.


\bibliographystyle{plain}
\bibliography{references,references1}

\end{document}